\documentclass[preprint,3p]{elsarticle}

\usepackage{amssymb}
\usepackage{amsmath}
\usepackage{empheq}
\usepackage{float}
\usepackage{physics}
\usepackage{siunitx}
\usepackage{subfig}
\usepackage{subfloat}
\usepackage{tikz}
\usetikzlibrary{arrows,positioning,shapes.geometric}
\DeclareSIUnit\Molar{M}

\DeclareMathOperator*{\argmin}{arg\,min\,}
\journal{Computer Methods in Applied Mechanics and Engineering}

\begin{document}

\begin{frontmatter}



\title{Influence of cellular mechano-calcium feedback in numerical models of cardiac electromechanics}


\author[a]{Irena Radišić\corref{cor1}}
\ead{irena.radisic@polimi.it}
\author[a]{Francesco Regazzoni}
\ead{francesco.regazzoni@polimi.it}
\author[a]{Michele Bucelli}
\ead{michele.bucelli@polimi.it}
\author[a]{Stefano Pagani}
\ead{stefano.pagani@polimi.it}
\author[a]{Luca Dede'}
\ead{luca.dede@polimi.it}
\author[a,b]{Alfio Quarteroni}
\ead{alfio.quarteroni@polimi.it}

\cortext[cor1]{Corresponding author}
\affiliation[a]{organization={MOX -- Department of Mathematics, Politecnico di Milano},%Department and Organization
            addressline={Piazza Leonardo da Vinci, 32}, 
            city={Milano},
            postcode={20133},
            country={Italy}}
\affiliation[b]{organization={Mathematics Institute (Professor Emeritus), École Polytechnique Fédérale de Lausanne},%Department and Organization
            addressline={Av. Piccard}, 
            city={Lausanne},
            postcode={CH-1015},
            country={Switzerland}}

%% Abstract
\begin{abstract}
Multiphysics and multiscale mathematical models enable the non-invasive study of cardiac function. These models often rely on simplifying assumptions that neglect certain biophysical processes to balance fidelity and computational cost. In this work, we propose an eikonal-based framework that incorporates mechano-calcium feedback -- the effect of mechanical deformation on calcium-troponin buffering -- while introducing only negligible computational overhead. To assess the impact of mechano-calcium feedback at the organ level, we develop a bidirectionally coupled cellular electromechanical model and integrate it into two cardiac multiscale frameworks: a monodomain-driven model that accounts for geometric feedback on electrophysiology and the proposed eikonal-based approach, which instead neglects geometric feedback. By ensuring consistent cellular model calibration across all scenarios, we isolate the role of mechano-calcium feedback and systematically compare its effects against models without it. Our results indicate that, under baseline conditions, mechano-calcium feedback has minimal influence on overall cardiac function. However, its effects become more pronounced in altered force generation scenarios, such as inotropic modulation. Furthermore, we demonstrate that the eikonal-based framework, despite omitting other types of mechano-electric feedback, effectively captures the role of mechano-calcium feedback at significantly lower computational costs than the monodomain-driven model, reinforcing its utility in computational cardiology.
\end{abstract}


%% Keywords

\end{frontmatter}



\section{Introduction}
\label{sec:introduction}
Mathematical models of cardiac electromechanics have increasingly gained biophysical detail making them a trustworthy tool in computational cardiology \cite{trayanova_computational_2012,taylor_computational_2013,zahid_patient-derived_2016,arevalo_arrhythmia_2016,passini_human_2017,gray_patient-specific_2018,prakosa_personalized_2018,boyle_computationally_2019,niederer_computational_2019,frontera_slow_2022}, but their computational cost makes their use prohibitive if access to computing resources is limited \cite{niederer_simulating_2011}. Balancing the model's biophysical accuracy with its computational efficiency becomes an increasingly challenging task as more detail is added to the mathematical model, calling for innovative formulations of the fully coupled problems \cite{neic_efficient_2017,regazzoni_machine_2022}.
\par
Mathematical models of cardiac electromechanics aim at simulating the cardiac function by the coupling of both electrical and mechanical processes, namely electrophysiology and mechanics. Key mechanisms underlying cardiac electromechanical function are the \textit{excitation-contraction coupling} \cite{katz_physiology_2011}, which describes the onset of a substantial contractile force due to a sudden rise of intracellular calcium following an electrical signal, and the \textit{mechano-electric feedbacks} \cite{franz_mechano-electrical_1996,timmermann_integrative_2017}, a wide array of phenomena which describe the effect of mechanical alterations on cardiac electrophysiology and which under normal regimes provide a self-regulating mechanism for cardiac contraction, but whose dysfunction is implicated in arrhythmogenic effects \cite{kohl_cardiac_2003, ravelli_mechano-electric_2003,timmermann_arrhythmogenic_2019,salvador_role_2022}. While computer simulation has helped elucidate the role of the mechano-electric feedbacks on the cardiac function \cite{salvador_role_2022,zappon_integrated_2024,gerach_differential_2024}, it is challenging to incorporate a large number of them consistently while understanding their differential effects.
\par
{In this study we specifically focus on the impact of the \textit{mechano-calcium feedback} (MCF) on the cardiac function. The MCF is a type of mechano-electric feedback which describes the effect of the sarcomere stretch on the calcium-troponin binding affinity \cite{allen_calcium_1988,terkeurs_arrhythmogenic_2006} and, consequently, on calcium dynamics. This feedback has often been neglected in cardiac electromechanics modeling \cite{fedele_comprehensive_2023,salvador_role_2022}, although recent works have incorporated the MCF in single-cell simulations \cite{margara_-silico_2021,mazhar_detailed_2024}, providing physiological results for single-cell simulations and highlighting the importance of the feedback on inotropic risk assessment. Notwithstanding that the MCF is a cellular-level mechanism, a more realistic representation of its effects can only be ascertained through macroscopic tissue simulations, where dynamic and spatially heterogeneous tissue mechanical properties are present. {There have been studies on the qualitative effects of the inclusion of the MCF in multiscale simulations \cite{gerach_differential_2024}, however they have mainly focused on its consequences in terms of electrophysiological substrate heterogeneity.}}
\par
The principal drawback of using electromechanical models with respect to electrophysiology only, besides a more difficult interpretation of cause-effect relationships, is its high computational complexity that increases with the number mechano-electric feedback mechanisms that are included \cite{gerbi_monolithic_2018,quarteroni_mathematical_2019}. To deal with this complexity, suitable segregated-staggered solution schemes have been devised \cite{salvador_intergrid_2020,fedele_comprehensive_2023}, along with stabilization and interpolation techniques capable of dealing with sources of instability and loss of accuracy \cite{regazzoni_oscillation-free_2021,bucelli_preserving_2023,bucelli_robust_2024}. While these techniques make simulation feasible, their associated cost is non-negligible, and the incorporation of feedbacks remains computationally demanding. On the other hand, model formulations associated with lower computational costs needed to compute the numerical solutions, while possibly disregarding certain biophysical details, such as formulations based on the eikonal equation, remain a valid alternative \cite{neic_efficient_2017,stella_fast_2022}. However, it is important to understand the effects of what is being neglected, and, given a specific application in mind, to show that it has adequately small effects on the quantities of interest.
\par
{In this work we study the impact of the inclusion of the mechano-calcium feedback in multiscale monodomain-driven cardiac electromechanics simulations by constructing bidirectionally coupled cellular electromechanical models, and for the first time, to the best of the our knowledge, we compare systematically their effects against models that do not account for the feedback starting from the same cellular baseline conditions. Moreover, we propose an alternative multiscale electromechanical framework based on the eikonal equation incorporating the mechano-calcium feedback. We show the eikonal-based framework to be a valid computationally efficient alternative even for biophysically detailed simulations, as exemplified by simulations on a realistic left ventricle geometry, both in baseline conditions and under variations of the cellular mechanical parameters.}
\par
This work is structured as follows. In Section \ref{sec:models_methods} we give an overview of mathematical models of cardiac electromechanics, we derive both single-cell and multiscale electromechanical models capable of capturing the MCF, and discuss their numerical implementation. In Section \ref{sec:numerical_results} we show the numerical results for the derived models. Finally, in Section \ref{sec:discussion}, we summarize the main contributions of this work and the implications for future cardiac electromechanics modeling and simulation.

\section{Models and methods}
\label{sec:models_methods}
In this section we overview mathematical models of cardiac electromechanics, highlighting the modeling assumptions where relevant, and present the modeling contributions of this work. In Section \ref{sec_model_derivation} we derive a bidirectionally coupled and biophysically consistent cellular electromechanical model, capturing both the excitation-contraction coupling and the mechano-calcium feedback, and we generalize this procedure to a wider class of single calcium-buffer models. In Section \ref{sec:fully_coupled_multiscale} we present the fully coupled multiscale model of cardiac electromechanics employing the introduced cellular electromechanics model. In Section \ref{sec:eikonal_driven} we present a computationally efficient multiscale framework for eikonal-driven cardiac electromechanics incorporating the mechano-calcium feedback. Finally, in Section \ref{sec:numerical_approximation}, we propose a numerical scheme appropriate for the simulation of all the models previously introduced in the work.
\subsection{Single-cell model derivation}
\label{sec_model_derivation}
Let us consider a general cellular ionic model, represented by the following system of ordinary differential equations (ODEs):
\begin{subequations}
    \label{eq:ionic}
    \begin{empheq}[left=\empheqlbrace]{align}
    & \dv{\boldsymbol{w}}{t} = \boldsymbol{g}\left(u, \boldsymbol{w}\right), \hspace{0.75em} &&\mbox{in} \ (0,T], \label{eq:ionic_variables}
    \\
    & C_\mathrm{m}\dv{u}{t} + \mathcal{I}_\mathrm{ion}\left(u,\boldsymbol{w}\right) = \mathcal{I}_{\mathrm{app}}(t), \hspace{0.75em} &&\mbox{in} \ (0,T], \label{eq:ionic_cable}
    \end{empheq}
\end{subequations}
with suitable initial conditions, where $\boldsymbol{w}=(\boldsymbol{w}_\mathrm{G},\boldsymbol{w}_\mathrm{I})$ are the ionic variables consisting of the gating variables and ion concentrations, $u$ is the transmembrane potential,
\begin{equation*}
    \boldsymbol{w}_\mathrm{G}:[0,T]\mapsto \mathbb{R}^{n_{\boldsymbol{w}_\mathrm{G}}}, \quad \boldsymbol{w}_\mathrm{I}:[0,T]\mapsto \mathbb{R}^{n_{\boldsymbol{w}_\mathrm{I}}}, \quad u:[0,T]\mapsto \mathbb{R},
\end{equation*}
$C_\mathrm{m}$ is the membrane capacitance, and $\mathcal{I}_\mathrm{app}$ is an externally applied, time dependent current. We will consider ionic models that include among their unknowns the free cytosolic intracellular calcium concentration $w_\mathrm{Ca}=[\mathrm{Ca}^{2+}]_\mathrm{i}$, as it is the main driver of the excitation-contraction coupling \cite{ten_tusscher_alternans_2006,tomek_development_2019}.
\par
Regarding the force generation, we will consider models that accurately describe the subcellular sarcomere dynamics. Since sarcomere contraction is regulated by calcium binding to troponin located on the sarcomeres' regulatory units, a biophysically accurate model of active force generation must describe the regulatory units' kinetics. More specifically, we will rely on the mean-field model proposed by Regazzoni et al. in \cite{regazzoni_biophysically_2020}, henceforth denoted by RDQ20, although this procedure may be extended to other models which provide the calcium-troponin buffering rate, e.g. see \cite{land_model_2017,mazhar_detailed_2024,margara_-silico_2021}. The force generation model takes the form:
\begin{equation}
    \label{eq:force_sarcomere_state}
    \dv{\boldsymbol{y}}{t} = \boldsymbol{h}\left(\boldsymbol{y},w_\mathrm{Ca},{\mathrm{SL}},{\dv{{\mathrm{SL}}}{t}}\right), \hspace{0.75em} \mbox{in} \ (0,T],
\end{equation}
accompanied by appropriate initial conditions, where the rate of evolution of the state of contractile and regulatory proteins, represented by the state vector $\boldsymbol{y}:[0,T]\mapsto \mathbb{R}^{n_{\boldsymbol{y}}}$, depends explicitly on the intracellular calcium concentration $w_\mathrm{Ca}$, the sarcomere length $\mathrm{SL}$ and its derivative $\dv{\mathrm{SL}}{t}$. The RDQ20 model captures the length dependence of the calcium-troponin binding affinity through $\mathrm{SL}$ and $\dv{\mathrm{SL}}{t}$, thereby affecting the regulatory units' dynamics. The resulting active force depends directly on the sarcomere state variables and sarcomere length, and need not be included as a state variable of the sarcomere dynamics model,
\begin{equation}
    \label{eq:active_force}
    T_\mathrm{a} = T_\mathrm{a}\left(\boldsymbol{y},\mathrm{SL}\right).
\end{equation}
From Equations \eqref{eq:ionic} and \eqref{eq:force_sarcomere_state}, it is evident that the forward coupling of $w_\mathrm{Ca}$ from the ionic model \eqref{eq:ionic} into the sarcomere dynamics model $\eqref{eq:force_sarcomere_state}$ describes the excitation-contraction mechanism. It is the form commonly used in electromechanical models which neglect the MCF \cite{fedele_comprehensive_2023}. We aim at extending this formulation to also include MCF. In order to clarify this aspect, it is necessary to explore the internal kinetics of the ionic and force generation models, and to this end in the next section we address the calcium kinetics in the ionic model.
\subsubsection{Subcellular calcium kinetics}
\label{sec:subcell_ca_kinetics}
For ionic models represented as systems of ODEs, such as \eqref{eq:ionic}, ion concentration evolutions are governed by extracellular ionic fluxes mediated by ionic channels, and intracellular ionic fluxes through various cellular compartments, mimicking the spatial disposition of a cardiomyocyte \cite{ten_tusscher_model_2004,ohara_simulation_2011}. The sarcomeres, although not explicitly represented, are typically situated in the bulk cytosolic compartment \cite{ten_tusscher_alternans_2006,tomek_development_2019}, where they appear simply as one of the calcium buffers (troponin).

The chemical reaction of calcium binding to and unbinding from its buffer in order to form the calcium-buffer complex $\mathrm{Buf_c}$ is expressed by the stoichiometric relation:
\begin{equation}
\label{eq:chemical_buffer}
    \mathrm{Ca}^{2+}+\mathrm{Buf} \xrightleftharpoons[k_-]{k_+} \mathrm{Buf_c},
\end{equation}
where the forward and backward reaction rate constants $k_+$ and $k_-$ encode the dependence of the reaction rates on the temperature, the polarity, the geometry of the interacting substances as well as other factors \cite{keener_mathematical_2009}.
The law of mass action \cite{keener_mathematical_2009} yields the following relations for the reactants' concentrations:
\begin{subequations}
\label{eq:rates_buffer}
    \begin{equation}
    \label{eq:rates_buffer_1}
        \dv{[\mathrm{Buf_c}]}{t} = k_+[\mathrm{Ca}^{2+}][\mathrm{Buf}] - k_-[\mathrm{Buf_c}]
            \end{equation}
            \begin{equation}
            \label{eq:rates_buffer_2}
    \dv{[\mathrm{Buf}]}{t} = k_-[\mathrm{Buf_c}] - k_+[\mathrm{Ca}^{2+}][\mathrm{Buf}],
    \end{equation}
\end{subequations}
here we have used $[\ \cdot\ ]$ to denote the molar concentrations. Since the right-hand sides of \eqref{eq:rates_buffer_1} and \eqref{eq:rates_buffer_2} sum to zero, the quantity $[\mathrm{Buf_c}] + [\mathrm{Buf}]$ is conserved in time:
\begin{equation}
\label{eq:total_buffer}
    [\mathrm{Buf_c}] + [\mathrm{Buf}] = \mathrm{Buf}_{\mathrm{c,max}},
\end{equation}
and it corresponds to the maximum allowed cellular calcium-buffer complex concentration $\mathrm{Buf}_{\mathrm{c,max}}$. Since calcium buffers are macromolecules and are not transported across the cellular membrane, the quantity $\mathrm{Buf}_{\mathrm{c,max}}$ is bound by the macromolecules' concentration and is thus a physical parameter of the model dictated by the cellular structure. Following a quasi-static approximation for the calcium-buffer complex concentration, the expression \eqref{eq:rates_buffer_1} becomes:
\begin{equation}
    \label{eq:rates_buffer_eq_1}
    k_+[\mathrm{Ca}^{2+}][\mathrm{Buf}] = k_-[\mathrm{Buf_c}],
\end{equation}
and by combining it with expression \eqref{eq:total_buffer}, we get the concentration of the calcium-buffer complex in dependence of the calcium concentration, expressed by the Michaelis-Menten equation:
\begin{equation}
    \label{eq:michaelis_menten_buffers}
    [\mathrm{Buf_c}] = \frac{[\mathrm{Ca}^{2+}]\mathrm{Buf}_{\mathrm{c,max}}}{[\mathrm{Ca}^{2+}] + K_{\mathrm{Buf_c}}},
\end{equation}
where the constant $K_{\mathrm{Buf_c}} = \frac{k_-}{k_+}$, is the half-saturation constant for the buffer $\mathrm{Buf}$.
\par
The quasi-static assumption yielding \eqref{eq:michaelis_menten_buffers} is also called the instantaneous equilibrium assumption, and it holds if the process of calcium binding to and from its buffer is much faster than the time scale on which the calcium concentration changes due to calcium influx into the cytosol \cite{keener_mathematical_2009}. For a simple chemical reaction with few intermediate reaction steps this is a reasonable assumption. In reality, there are effects due to spatial concentration gradients and the proximity of the calcium-buffer binding sites from calcium influx and efflux sites \cite{colman_multi-scale_2022, keener_mathematical_2009}. This true spatial gradient cannot be present in a zero-dimensional ionic model, as the ones of type \eqref{eq:ionic_variables} we are considering. Some spatial resolution may be recovered by introducing various cellular compartments. However, in the individual cytosolic compartments the calcium concentration is uniform and thus no spatial gradient is present. A further, and perhaps more striking, flaw of the considered ionic models is the restrictions made on the reaction rate constants $k_+$ and $k_-$. As mentioned above, the rates depend on the temperature but also the geometry of the interacting molecules, which becomes non-negligible when considering the large strain variation and subsequent microscopic stretch associated with cardiac contraction \cite{dobesh_cooperative_2002,rice_approaches_2004}. Notwithstanding these considerations, the reaction rates in the ionic models we will be considering are constant and fixed.

\subsubsection{Feedback condition}
\label{sec:feedback_condition}
As discussed in Section \ref{sec:subcell_ca_kinetics}, ionic models typically use a Michaelis-Menten approximation with constant reaction rates for the calcium buffering kinetics \cite{luo_dynamic_1994,ten_tusscher_model_2004,ohara_simulation_2011}, including the calcium-troponin buffering reaction. However, as we will show, the calcium-troponin reaction is also described by the active force generation model as part of the regulatory kinetics, possibly creating an inconsistency between the two microscopic descriptions. In order to formulate biophysically consistent models, these discrepancies can be solved by considering appropriate coupling conditions.
\par
To this end, let us consider the RDQ20 activation model. The fraction of regulatory troponin units bound to calcium depends on the sarcomere length and is given by:
\begin{equation}
\label{eq:bound_expr}
    \mathcal{B}\left(\boldsymbol{y},\mathrm{SL}\right) = \chi_\mathrm{so}(\mathrm{SL})\mathcal{B}_\mathrm{so}(\boldsymbol{y}) + \left(1-\chi_\mathrm{so}(\mathrm{SL})\right)\mathcal{B}_\mathrm{nso}(\boldsymbol{y}),
\end{equation}
where $\mathcal{B}_\mathrm{so}(\boldsymbol{y})$ and $\mathcal{B}_\mathrm{nso}(\boldsymbol{y})$ are the fractions of bound troponin units in the single-overlap and non-single-overlap zones,
\begin{equation*}
    \mathcal{B}_\mathrm{so}:\mathbb{R}^{n_{\boldsymbol{y}}}\mapsto [0,1], \quad \mathcal{B}_\mathrm{nso}:\mathbb{R}^{n_{\boldsymbol{y}}}\mapsto [0,1],
\end{equation*}
and $\chi_\mathrm{so}(\mathrm{SL})$ is the fraction of the sarcomere in the single-overlap zone, and has the following expression:
\begin{equation}
\label{eq:def_chiso}
\chi_{\mathrm{so}}(\mathrm{SL})= \begin{dcases}0, & \mbox { if } \mathrm{SL} \leq \mathrm{L_A}, \\ \frac{2\left(\mathrm{SL}-\mathrm{L_A}\right)}{\mathrm{L_M}-\mathrm{L_H}}, & \mbox { if } \mathrm{L_A}<\mathrm{SL} \leq \mathrm{L_M}, \\ \frac{\mathrm{SL}+\mathrm{L_M}-2 \mathrm{L_A}}{\mathrm{L_M}-\mathrm{L_H}}, & \mbox { if } \mathrm{L_M}<\mathrm{SL} \leq 2 \mathrm{L_A}-\mathrm{L_H}, \\ 1, & \mbox { if } 2 \mathrm{L_A}-\mathrm{L_H}<\mathrm{SL} \leq 2 \mathrm{L_A}+\mathrm{L_H}, \\ \frac{\mathrm{L_M}+2 \mathrm{L_A}-\mathrm{SL}}{\mathrm{L_M}-\mathrm{L_H}}, & \mbox { if } 2 \mathrm{L_A}+\mathrm{L_H}<\mathrm{SL} \leq 2 \mathrm{L_A}+\mathrm{L_M}, \\ 0, & \mbox { if } \mathrm{SL}>2 \mathrm{L_A}+\mathrm{L_M}.\end{dcases}
\end{equation}
Recalling that the calcium-bound troponin concentration is given by the fraction of bound regulatory troponin units multiplied by the total concentration of the troponin units, which coincides with the maximal possible calcium-troponin concentration ${\mathrm{Tn}}_{\mathrm{c,max}}$, it follows:
\begin{equation}
    [\mathrm{Ca}^{2+}]_\mathrm{Tn}\left(\boldsymbol{y},\mathrm{SL}\right) = {\mathrm{Tn}}_{\mathrm{c,max}}\mathcal{B}\left(\boldsymbol{y},\mathrm{SL}\right).
\end{equation}

Deriving in time and invoking \eqref{eq:bound_expr}, one gets,
\begin{equation}
\label{eq:derived_ca_tn_active_force}
\begin{aligned}
    \dv{[\mathrm{Ca}^{2+}]_\mathrm{Tn}}{t}={\mathrm{Tn}}_{\mathrm{c,max}}\biggl(&\chi_\mathrm{so}(\mathrm{SL})\dv{\mathcal{B}_\mathrm{so}(\boldsymbol{y})}{t} + \left(1-\chi_\mathrm{so}(\mathrm{SL})\right)\dv{\mathcal{B}_\mathrm{nso}(\boldsymbol{y})}{t}\\ &+ \left( \mathcal{B}_\mathrm{so}\left(\boldsymbol{y}) - \mathcal{B}_\mathrm{nso}(\boldsymbol{y}\right)\right) \dv{\mathcal{B}}{\mathrm{SL}}\dv{\mathrm{SL}}{t} \biggr),
\end{aligned}
\end{equation}
with $\dv{\mathcal{B}}{\mathrm{SL}}$ being a piecewise constant function of $\mathrm{SL}$. Considering the transition rates from a bound to an unbound state and vice-versa for the regulatory units from the RDQ20 model:
\begin{equation}
    \label{eq:transition_bound}
    \dv{\mathcal{B}_\star}{t} = \frac{k_\mathrm{off}}{k_\mathrm{d}(\mathrm{SL})}[\mathrm{Ca}^{2+}]_\mathrm{i}\left(1-\mathcal{B}_\star\right) - k_\mathrm{off}\mathcal{B}_\star, \qquad \star\in\left\{\mathrm{so},\mathrm{nso}\right\},
\end{equation}
expression \eqref{eq:derived_ca_tn_active_force} becomes:
\begin{equation}
\label{eq:derived_ca_tn_active_force_detail}
\begin{aligned}
    \dv{[\mathrm{Ca}^{2+}]_\mathrm{Tn}}{t}={\mathrm{Tn}}_{\mathrm{c,max}}\biggl(&\chi_\mathrm{so}(\mathrm{SL}) \left(\frac{k_\mathrm{off}}{k_\mathrm{d}(\mathrm{SL})}[\mathrm{Ca}^{2+}]_\mathrm{i}\left(1-\mathcal{B}_\mathrm{so}(\boldsymbol{y})\right) - k_\mathrm{off}\mathcal{B}_\mathrm{so}(\boldsymbol{y})\right) \\ & +  \left(1-\chi_\mathrm{so}(\mathrm{SL})\right)\left(\frac{k_\mathrm{off}}{k_\mathrm{d}(\mathrm{SL})}[\mathrm{Ca}^{2+}]_\mathrm{i}\left(1-\mathcal{B}_\mathrm{nso}(\boldsymbol{y})\right) - k_\mathrm{off}\mathcal{B}_\mathrm{nso}(\boldsymbol{y})\right)\\ & + \left( \mathcal{B}_\mathrm{so}\left(\boldsymbol{y}) - \mathcal{B}_\mathrm{nso}(\boldsymbol{y}\right)\right) \dv{\mathcal{B}}{\mathrm{SL}}\dv{\mathrm{SL}}{t} \biggr),
\end{aligned}
\end{equation}
thus allowing to express the calcium-troponin buffering rate as a function of the sarcomere model's state variables, the sarcomere length and its derivative:
\begin{equation}
\label{eq:catn_dt_active_force}
    \dv{[\mathrm{Ca}^{2+}]_\mathrm{Tn}}{t} = y_{\mathrm{dCa_{Tn}}}\left(\boldsymbol{y},\mathrm{SL},\dv{\mathrm{SL}}{t}\right).
\end{equation}
In order to find the same quantity in a given ionic model, we start from the total bulk cytosolic calcium concentration, noting that it may either be free, denoted with $[\mathrm{Ca}^{2+}]_\mathrm{i}$, bound to troponin $[\mathrm{Ca}^{2+}]_\mathrm{Tn}$, or bound to any other generic cytosolic buffer $[\mathrm{Ca}^{2+}]_\mathrm{Buf}$, yielding:
\begin{equation}
\label{eq:total_calcium}
    [\mathrm{Ca}^{2+}]_\mathrm{tot} = [\mathrm{Ca}^{2+}]_\mathrm{i} + [\mathrm{Ca}^{2+}]_\mathrm{Tn} + [\mathrm{Ca}^{2+}]_\mathrm{Buf},
\end{equation}
which, by differentiating in time, becomes:
\begin{equation}
\label{eq:total_calcium_time_derivative}
    \dv{[\mathrm{Ca}^{2+}]_\mathrm{tot}}{t} = \dv{[\mathrm{Ca}^{2+}]_\mathrm{i}}{t} + \dv{[\mathrm{Ca}^{2+}]_\mathrm{Tn}}{t} + \dv{[\mathrm{Ca}^{2+}]_\mathrm{Buf}}{t}.
\end{equation}
Previous works \cite{ten_tusscher_model_2004,ohara_simulation_2011,luo_dynamic_1994} apply relation \eqref{eq:michaelis_menten_buffers} to express the rate of change of both the buffer-bound and troponin-bound calcium. Instead, we only do so for the buffer-bound calcium, and instead use \eqref{eq:catn_dt_active_force} for the troponin-bound calcium, thereby relaxing its instantaneous equilibrium assumption. Rearranging the terms in \eqref{eq:total_calcium_time_derivative} and applying the above considerations, it becomes:
\begin{equation}
\label{eq:free_intracellular_pre}
    \dv{[\mathrm{Ca}^{2+}]_\mathrm{i}}{t} = \dv{[\mathrm{Ca}^{2+}]_\mathrm{tot}}{t} - y_{\mathrm{dCa_{Tn}}}\left(\boldsymbol{y},\mathrm{SL},\dv{\mathrm{SL}}{t}\right)-\dv{}{t}\left(\frac{[\mathrm{Ca}^{2+}]_\mathrm{i}\mathrm{Buf}_{\mathrm{c,max}}}{[\mathrm{Ca}^{2+}]_\mathrm{i} + K_{\mathrm{Buf_c}}}\right).
\end{equation}
Mass conservation dictates that that the total calcium variation can only be due to transmembrane or transcompartmental calcium-carrying currents, which may be formulated as:
\begin{equation}
    \label{eq:mass_conservation_calcium}
    \dv{[\mathrm{Ca}^{2+}]_\mathrm{tot}}{t} = \sum_{i=1}^{N_\mathrm{curr}}I_{\mathrm{Ca},i},
\end{equation}
where $\left\{I_{\mathrm{Ca},i}\right\}_{i = 1}^{N_\mathrm{curr}}$ is the set of all calcium-carrying currents of the bulk cytosol, depending on the ionic model and expressed either in current or in flux form. By substituting \eqref{eq:mass_conservation_calcium} into \eqref{eq:free_intracellular_pre}, applying the chain rule, and rearranging the terms in \eqref{eq:free_intracellular_pre}, we obtain a modified equation for the evolution of the free intracellular cytosolic calcium, which reads:
\begin{equation}
\label{eq:dcai_modified}
\begin{aligned}
    \dv{[\mathrm{Ca}^{2+}]_\mathrm{i}}{t} = &\left({ 1 + \frac{\mathrm{Buf}_{\mathrm{c,max}}K_{\mathrm{Buf_c}}}{\left([\mathrm{Ca}^{2+}]_\mathrm{i} + K_{\mathrm{Buf_c}}\right)^2}}\right)^{-1}\left({\sum_{i=1}^{N_\mathrm{curr}}I_{\mathrm{Ca},i} - y_{\mathrm{dCa_{Tn}}}\left(\boldsymbol{y},\mathrm{SL},\dv{\mathrm{SL}}{t}\right)}\right),
\end{aligned}
\end{equation}
as opposed to the evolution law obtained by assuming instantaneous equilibrium \eqref{eq:michaelis_menten_buffers} for the calcium-troponin buffering:
\begin{equation}
\label{eq:dcai_original}
\begin{aligned}
    \dv{[\mathrm{Ca}^{2+}]_\mathrm{i}}{t} = &\left({ 1 + \frac{\mathrm{Tn}_{\mathrm{c,max}}K_{\mathrm{Tn_c}}}{\left([\mathrm{Ca}^{2+}]_\mathrm{i} + K_{\mathrm{Tn_c}}\right)^2} + \frac{\mathrm{Buf}_{\mathrm{c,max}}K_{\mathrm{Buf_c}}}{\left([\mathrm{Ca}^{2+}]_\mathrm{i} + K_{\mathrm{Buf_c}}\right)^2}}\right)^{-1}{\sum_{i=1}^{N_\mathrm{curr}}I_{\mathrm{Ca},i}}.
\end{aligned}
\end{equation}
Equation \eqref{eq:dcai_modified} depends on $\boldsymbol{y},\mathrm{SL},\pdv{\mathrm{SL}}{t}$ through the term $y_{\mathrm{dCa_{Tn}}}$, as opposed to Equation \eqref{eq:dcai_original} which depends only on the ionic variables $\boldsymbol{w}$. Thus we obtain a fully coupled cellular electromechanical model:
\begin{subequations}
    \label{eq:full_em_cell}
    \begin{empheq}[left=\empheqlbrace]{align}
    & \dv{\boldsymbol{w}}{t} = \tilde{\boldsymbol{g}}\left(u, \boldsymbol{w}, y_{\mathrm{dCa_{Tn}}}\left(\boldsymbol{y},{\mathrm{SL}},{\dv{{\mathrm{SL}}}{t}}\right)\right), \label{eq:ionic_variables_mod}
    \\
    & \dv{\boldsymbol{y}}{t} = \boldsymbol{h}\left(\boldsymbol{y},w_\mathrm{Ca},{\mathrm{SL}},{\dv{{\mathrm{SL}}}{t}}\right), \label{eq:sarcomere_variables_mod}
    \\
    & C_\mathrm{m}\dv{u}{t} + \mathcal{I}_\mathrm{ion}\left(u,\boldsymbol{w}\right) = \mathcal{I}_{\mathrm{app}}(t),\label{eq:ionic_cable_mod}
    \end{empheq}
\end{subequations}
where $\tilde{\boldsymbol{g}}$ is obtained by replacing Equation \eqref{eq:dcai_original} in \eqref{eq:ionic} with Equation \eqref{eq:dcai_modified}.
\par
We conclude this section by highlighting how expressions \eqref{eq:transition_bound} and \eqref{eq:derived_ca_tn_active_force_detail} give an insight into the instantaneous equilibrium assumptions discussed in Section \ref{sec:subcell_ca_kinetics}. Indeed, from the RDQ20 model \cite{regazzoni_biophysically_2020}, $k_\mathrm{off}$ corresponds, with appropriate rescaling, to the forward and backward reaction rates $k_+$ and $k_-$ presented in Equation \eqref{eq:chemical_buffer} for the calcium-troponin binding reaction. Thus, for $k_\mathrm{off}\rightarrow\infty$ the reaction, for a simple single-step approximation, becomes in instantaneous equilibrium, satisfying the initial assumptions of the ionic model. However, the half-saturation constant $K_{\mathrm{Tn_c}}$, which now depends on the kinetics \eqref{eq:transition_bound}, depends on the sarcomere length through the dependence of the dissociation constant $k_\mathrm{d}(\mathrm{SL})$ on the sarcomere length $\mathrm{SL}$, shifting the equilibrium point. More specifically, in the RDQ20 model the calcium-troponin binding is not necessarily single-step \cite{regazzoni_biophysically_2020}, but instead it has an intermediate step given by the permissivity transition whose rate is scaled by $k_\mathrm{basic}$, making it rate-limiting for $k_\mathrm{off}\rightarrow\infty$. Thus, more accurately, the instantaneous equilibrium limit is reached for $k_\mathrm{off},k_\mathrm{basic}\rightarrow\infty$ with appropriate scaling. Regardless of whether the equilibrium point is reached instantaneously or not, it is still shifted due to the effect of the sarcomere length.

\subsubsection{Single-buffer models}
\label{sec:single_buffer}
Some ionic models do not provide separate kinetic descriptions for different bulk cytosolic buffers, as for example the ten Tusscher-Panfilov (TTP06) model \cite{ten_tusscher_alternans_2006}. In these cases there is a single aggregate buffer population $\overline{\mathrm{Buf}}$ characterized by a single half-saturation constant $K_{\overline{\mathrm{Buf}}_\mathrm{c}}$ and a single maximal buffer concentration $\overline{\mathrm{Buf}}_{\mathrm{c,max}}$. Still, it is attractive from a computational point of view to implement this kind of coupling for such models, as oftentimes they provide less expensive but physically meaningful alternatives to more complete ionic models \cite{tomek_development_2019,fedele_comprehensive_2023} employed in multiscale simulations. In order to implement the feedback condition from the sarcomere model, it is necessary to separate the troponin from the other types of buffers in the ionic model. To this end, we split the buffer maximal aggregate calcium-buffer complex concentration into its troponin and other parts as:
\begin{equation}
    \label{eq:split_buffer}
    \overline{\mathrm{Buf}}_{\mathrm{c,max}} = {\mathrm{Buf}}_{\mathrm{c,max}} + {\mathrm{Tn}}_{\mathrm{c,max}},
\end{equation}
where ${\mathrm{Tn}}_{\mathrm{c,max}} >0$ is a free, positive model parameter representing the maximal calcium-troponin buffer concentration, to be subsequently calibrated. This split is trivially equivalent to the original ionic model under the assumptions that both buffer populations have the same reaction kinetics, which is not true in general \cite{smith_calcium_2019}. With these considerations the modified free intracellular cytosolic calcium concentration evolution equation, similarly to \eqref{eq:dcai_modified}, reads:
\begin{equation*}
    \dv{[\mathrm{Ca}^{2+}]_\mathrm{i}}{t} = \left({ 1 + \frac{\left( \overline{\mathrm{Buf}}_{\mathrm{c,max}} - {\mathrm{Tn}}_{\mathrm{c,max}} \right)K_{\mathrm{Buf_c}}}{\left([\mathrm{Ca}^{2+}]_\mathrm{i} + K_{\mathrm{Buf_c}}\right)^2}}\right)^{-1}\left({\sum_{i=1}^{N_\mathrm{curr}}I_{\mathrm{Ca},i} - y_{\mathrm{dCa_{Tn}}}\left(\boldsymbol{y},\mathrm{SL},\dv{\mathrm{SL}}{t}\right)}\right).
\end{equation*}
Reassuming, regarding the split \eqref{eq:split_buffer}, under the assumption that both buffer populations have the same reaction kinetics the split ionic model is equivalent to the unsplit one. Also, trivially, the model incorporating the feedback is equivalent to the original one for $\mathrm{Tn_{c,max}}=0$. In principle, the choice of the physical parameter ${\mathrm{Tn}}_{\mathrm{c,max}}$ could come from calcium buffer data and measurements for the human cardiomyocyte. However the lack thereof \cite{fabiato_calcium-induced_1983,gao_myofilament_1994,robertson_time-course_1981} makes its a priori choice difficult and calls for specific parameter calibration and inverse estimation techniques, as we will see in Section \ref{sec:cellular_electromechanical_model}.
\subsection{Fully coupled multiscale cardiac electromechanics}
\label{sec:fully_coupled_multiscale}
In this section we present a fully coupled multiscale model, henceforth denoted as by M+MCF, of the left ventricle electromechanics, starting from the ones presented in \cite{fedele_comprehensive_2023,regazzoni_cardiac_2022}, but including the mechano-calcium feedback effect expressed by \eqref{eq:catn_dt_active_force}. The full electromechanical model is cast in the unloaded conductive muscular domain $\Omega_0\subset\mathbb{R}^3$. We consider the domain $\Omega_0$ to be the volume occupied by the left ventricle with an artificial boundary obtained by neglecting its uppermost part, as in \cite{regazzoni_cardiac_2022}. We split the boundary $\partial\Omega_0$ into the endocardial surface $\Gamma_0^\mathrm{endo}$, the epicardial surface $\Gamma_0^\mathrm{epi}$, and the basal surface $\Gamma_0^\mathrm{base}$ corresponding to the artificial section. With this definition of the computational domain, the full electromechanical model reads:
\begin{subequations}
    \label{eq:full_em_space}
    \begin{empheq}[left=\empheqlbrace]{align}
    & J\chi_\mathrm{m}\left[C_\mathrm{m}\pdv{u}{t} + \mathcal{I}_\mathrm{ion}(u,\boldsymbol{w})\right]-\nabla\cdot\left(J\mathbf{F}^{-1}\mathbf{D}\mathbf{F}^{-T}\nabla u\right) = J\chi_\mathrm{m}\mathcal{I}_\mathrm{app}(t), \hspace{-0.5cm}& \mbox{in} \ \Omega_0\cross(0,T],
    \label{eq:monodomain}
    \\
    & \pdv{\boldsymbol{w}}{t} = \tilde{\boldsymbol{g}}\left(u, \boldsymbol{w}, y_{\mathrm{dCa_{Tn}}}\left(\boldsymbol{y},{\mathrm{SL}},{\pdv{{\mathrm{SL}}}{t}}\right)\right), & \mbox{in} \ \Omega_0\cross(0,T],
    \label{eq:ionic_variables_space}
    \\
    & \pdv{\boldsymbol{y}}{t} = \boldsymbol{h}\left(\boldsymbol{y},w_\mathrm{Ca},{\mathrm{SL}},{\pdv{{\mathrm{SL}}}{t}}\right), & \mbox{in} \ \Omega_0\cross(0,T], \label{eq:sarcomere_variables_space}
    \\
    & \rho_\mathrm{s}\pdv[2]{\mathbf{d}}{t} - \nabla\cdot\mathbf{P}\left(\mathbf{d},T_\mathrm{a}\left(\boldsymbol{y},\mathrm{SL}\right)\right) = 0, & \mbox{in} \ \Omega_0\cross(0,T],
    \label{eq:mechanics}
    \\
    & \mathcal{C}\left(p_\mathrm{LV},V_\mathrm{LV}(\mathbf{d}),t\right) = 0, & \mbox{in} \ \Omega_0\cross(0,T],
    \label{eq:circulation}
    \\
    & \left(J\mathbf{F}^{-1}\mathbf{D}\mathbf{F}^{-T}\nabla u\right)\cdot\mathbf{N}=0, & \mbox{on} \ \partial\Omega_0\cross(0,T],
    \label{eq:bc_monodomain}
    \\
    & \mathbf{P}\left(\mathbf{d},T_\mathrm{a}\left(\boldsymbol{y},\mathrm{SL}\right)\right)\mathbf{N} = -p_\mathrm{LV}(t)J\mathbf{F}^{-T}\mathbf{N}, & \mbox{on} \ \Gamma_0^\mathrm{endo}\cross(0,T], \label{eq:bc_endo}
    \\
    & \mathbf{P}\left(\mathbf{d},T_\mathrm{a}\left(\boldsymbol{y},\mathrm{SL}\right)\right)\mathbf{N} + \mathbf{K}^\mathrm{epi}\mathbf{d} + \mathbf{C}^\mathrm{epi}\pdv{\mathbf{d}}{t} = \mathbf{0}, & \mbox{on} \ \Gamma_0^\mathrm{epi}\cross(0,T],
    \label{eq:bc_epi}
    \\
    & \mathbf{P}\left(\mathbf{d},T_\mathrm{a}\left(\boldsymbol{y},\mathrm{SL}\right)\right)\mathbf{N} = p_\mathrm{LV}(t)||J\mathbf{F}^{-T}\mathbf{N}|| \frac{\int_{{\Gamma}^\mathrm{endo}_0}J\mathbf{F}^{-T}\mathbf{N} \mathrm{d}A }{\int_{\Gamma^\mathrm{base}_0}||J\mathbf{F}^{-T}\mathbf{N}||\mathrm{d}A}, & \mbox{on} \ \Gamma_0^\mathrm{base}\cross(0,T],
    \label{eq:bc_base}
    \\
    & \left(u,\boldsymbol{w},\boldsymbol{y},\mathbf{d},\dot{\mathbf{d}}\right)(0)=\left(u_0,\boldsymbol{w}_0,\boldsymbol{y}_0,\mathbf{d}_0,\mathbf{0}\right), & \mbox{in} \ \Omega_0,
    \label{eq:ic}
    \end{empheq}
\end{subequations}
where $\mathrm{SL}$ denotes the sarcomere length, computed as
\begin{equation*}
    \mathrm{SL} = \mathrm{SL}_0\norm{\mathbf{F}\mathbf{f}_0}.
\end{equation*}
The system \eqref{eq:full_em_space} is completed with suitable initial conditions \eqref{eq:ic}, and boundary conditions \eqref{eq:bc_endo},\eqref{eq:bc_epi},\eqref{eq:bc_base} for \eqref{eq:monodomain} and \eqref{eq:mechanics}. Condition \eqref{eq:bc_monodomain} is the no flux condition for the transmembrane potential, \eqref{eq:bc_endo} is the intraventricular pressure $p_\mathrm{LV}$ exerted on the endocardial surface and determined by a circulation model as in \cite{regazzoni_machine_2020}, \eqref{eq:bc_epi} models the effects of the pericardial sac \cite{pfaller_importance_2019}, and \eqref{eq:bc_base} is the energy-consistent boundary condition proposed in \cite{regazzoni_machine_2020}. Besides the ionic variables $\boldsymbol{w}$, the sarcomere state $\boldsymbol{y}$, and the transmembrane potential $u$, system \eqref{eq:full_em_space} describes the evolution of the macroscopic total Lagrangian tissue displacement $\mathbf{d}$, obtained from the momentum conservation equation \eqref{eq:mechanics} \cite{gurtin_mechanics_2010}. The myocardium is modeled as an orthotropic active material \cite{holzapfel_constitutive_2009, ambrosi_active_2012}, with its principal material directions $\left\{\mathbf{f}_0,\mathbf{s}_0,\mathbf{n}_0\right\}$, determined using the Bayer et al. rule-based algorithm \cite{bayer_novel_2012}, rotating transmurally in the myocardium \cite{piersanti_modeling_2021}. The myocardial passive mechanical properties are modeled as hyperelastic:
\begin{equation}
    \label{eq:stress_decomp}
\mathbf{P}\left(\mathbf{d},T_\mathrm{a}\left(\boldsymbol{y},\mathrm{SL}\right)\right) = \pdv{\mathcal{W}(\mathbf{F})}{\mathbf{F}} + T_\mathrm{a}\left(\boldsymbol{y},\mathrm{SL}\right)\frac{\mathbf{F}\mathbf{f}_0\otimes\mathbf{f}_0}{||\mathbf{F}\mathbf{f}_0||},
\end{equation}
where $\mathcal{W}$ is the strain energy density, for which we employ the exponential Usyk et al. model \cite{usyk_computational_2002}, and $T_\mathrm{a}$ is the active contribution of microscopic force generation mechanisms acting in the principle muscle fiber direction \cite{regazzoni_cardiac_2022}. Regarding the electrical propagation, Equation \eqref{eq:ionic_cable_mod} is replaced with the monodomain \cite{colli_franzone_mathematical_2014} equation \eqref{eq:monodomain} written in the reference domain with the conductivity tensor $\mathbf{D}$ and cell membrane surface to volume ratio $\chi_\mathrm{m}$ being modified by the deformation gradient $\mathbf{F}$ and its determinant $J$, thus incorporating the influence of the mechanics on the signal propagation. The myocardium is initially stimulated by a locally applied current $\mathcal{I}_\mathrm{app}(t)$.
\par
The mechanical problem is closed by coupling it with a circulation model \eqref{eq:circulation} which yields a pressure-volume relationship dependent of the phase of the cardiac cycle. We use, as in \cite{stella_fast_2022}, a Lagrange multiplier based isovolumic constraint for the two isovolumic phases and a two-element Windkessel afterload model for the ejection phase, which relates the pressure and the volume as:
\begin{subequations}
    \label{eq:windkessel}
    \begin{empheq}[left=\empheqlbrace]{align}
    & C_\mathrm{WK}\dv{p_\mathrm{LV}}{t}=-\frac{p_\mathrm{LV}}{R_\mathrm{WK}}-\dv{V_\mathrm{LV}}{t}, \hspace{0.75em} t\in\left(t_\mathrm{AVO},t_\mathrm{AVC}\right)\\
    &p_\mathrm{LV}(t_\mathrm{AVO})=p_\mathrm{AVO},
    \end{empheq}
\end{subequations}
where $C_\mathrm{WK},R_\mathrm{WK},p_\mathrm{AVO}$ are model parameters. For more details on the circulation model we refer to \cite{stella_fast_2022}. Despite its relative simplicity, the benefit of this afterload model is the shorter time to reach a limit cycle and easier calibration with respect to more detailed models \cite{fedele_comprehensive_2023}.
\par
Due to the different depolarization times owed to the spatial propagation of the electrical signal $u$ \eqref{eq:monodomain}, the ionic and sarcomere state variables $\boldsymbol{w}$ and $\boldsymbol{y}$ are space-dependent, and their heterogeneity is further amplified by the MCF in Equation \eqref{eq:ionic_variables_space} and its dependence on the local tissue displacement through the sarcomere length $\mathrm{SL}$. With respect to the single-cell model \eqref{eq:full_em_cell}, the multiscale model \eqref{eq:full_em_space} takes into account the direct feedback from the mechanics to the ionic model through the sarcomere length and its time derivative.
\par
The numerical approximation of system \eqref{eq:full_em_space} {with classical methods} is computationally expensive, due to the fast timescales involved and high spatial resolution required in order to capture the propagating wavefronts in \eqref{eq:monodomain}, along with the nonlinearities due to the material models and couplings involved \cite{colli_franzone_wavefront_1990, quarteroni_mathematical_2019}, limiting the use of \eqref{eq:full_em_space} with highly detailed ionic models appropriate for applications such as \cite{tomek_development_2019}. It is therefore of interest to reduce the computational costs associated with these types of problems.
\subsection{Eikonal-driven cardiac multiscale electromechanics}
\label{sec:eikonal_driven}
A less computationally demanding approach to model cardiac electromechanics consists of pre-computing the activation sequences of the myocardium by means of an Eikonal model \cite{colli_franzone_mathematical_2014, neic_efficient_2017, stella_fast_2022}. In particular, we use the Eikonal-Diffusion model \cite{colli_franzone_wavefront_1990} as in \cite{stella_fast_2022}, which consists of finding the activation time $\psi:\Omega_0\rightarrow\mathbb{R}$ for each point of the conductive domain $\Omega_0$ by solving:
\begin{subequations}
    \label{eq:eikonal}
    \begin{empheq}[left=\empheqlbrace]{align}
    & c_0\sqrt{\nabla\psi\cdot\frac{1}{\chi C_\mathrm{m}}\mathbf{D}\nabla\psi} - \nabla\cdot\left(\frac{1}{\chi C_\mathrm{m}}\mathbf{D}\nabla\psi\right) = 1, \hspace{0.75em} &&\mbox{in} \ \Omega_0,
    \label{eq:eik_system}
    \\
    & \left(\frac{1}{\chi C_\mathrm{m}}\mathbf{D}\nabla\psi\right)\cdot\mathbf{N} = 0, \hspace{0.75em} &&\mbox{on} \ \partial\Omega_0\setminus\partial\Omega_\mathrm{a}, \label{eq:eik_boundary}
    \\
    & \psi = \psi_\mathrm{a}, \hspace{0.75em} &&\mbox{on} \ \partial\Omega_\mathrm{a}.
    \label{eq:eik_dirichlet}
    \end{empheq}
\end{subequations}
where $c_0$ determines the wavefront velocity. In this manner the mechano-electric feedbacks on the depolarization related to the moving geometry are disregarded. It is still, however, possible to include the effect of the mechano-calcium feedback, by solving the modified system, henceforth denoted as E+MCF, given the solution $\psi$ of system \eqref{eq:eikonal}:
\begin{subequations}
    \label{eq:eik_em_space}
    \begin{empheq}[left=\empheqlbrace]{align}
    & \pdv{\boldsymbol{w}}{t} = \tilde{\boldsymbol{g}}\left(u, \boldsymbol{w}, y_{\mathrm{dCa_{Tn}}}\left(\boldsymbol{y},{\mathrm{SL}},{\dv{{\mathrm{SL}}}{t}}\right)\right), \hspace{0.75em} &&\mbox{in} \ \Omega_0\cross(0,T],
    \label{eq:ionic_variables_space_mcf}
    \\
    & \pdv{\boldsymbol{y}}{t} = \boldsymbol{h}\left(\boldsymbol{y},w_\mathrm{Ca},{\mathrm{SL}},{\pdv{{\mathrm{SL}}}{t}}\right), \hspace{0.75em} &&\mbox{in} \ \Omega_0\cross(0,T], \label{eq:sarcomere_variables_space_mcf}
    \\
    & C_\mathrm{m}\dv{u}{t} + \mathcal{I}_\mathrm{ion}\left(u,\boldsymbol{w}\right) = \mathcal{I}_{\mathrm{app}}(t-\psi), \hspace{0.75em} &&\mbox{in} \ \Omega_0\cross(0,T],
    \label{eq:ionic_cable_mcf}
    \\
    & \rho_\mathrm{s}\pdv[2]{\mathbf{d}}{t} - \nabla\cdot\mathbf{P}\left(\mathbf{d},T_\mathrm{a}\left(\boldsymbol{y},\mathrm{SL}\right)\right) = 0, \hspace{0.75em} &&\mbox{in} \ \Omega_0\cross(0,T],
    \label{eq:mechanics_mcf}
    \\
    & \mathcal{C}\left(p_\mathrm{LV},V_\mathrm{LV}(\mathbf{d}),t\right) = 0, \hspace{0.75em} &&\mbox{in} \ \Omega_0\cross(0,T],
    \label{eq:circulation_eikonal}
    \end{empheq}
\end{subequations}
with appropriate initial conditions and boundary conditions \eqref{eq:bc_endo}, \eqref{eq:bc_epi} and \eqref{eq:bc_base}. The depolarizing diffusive currents in the monodomain equation are surrogated in \eqref{eq:ionic_cable_mcf} by the locally applied current
\begin{equation}
    \label{eq:applied_current}
    \mathcal{I}_\mathrm{app}(t) = \Bar{\mathcal{I}}_\mathrm{app}\mathbb{I}_{\left\{0<t<t_\mathrm{app} \right\}}, \quad t\in\left[0,T \right],
\end{equation}
where $\Bar{\mathcal{I}}_\mathrm{app}$ is the current amplitude and $t_\mathrm{app}<T$ is its duration, whose application time is shifted by the depolarizing current's arrival time approximated by the computed activation time $\psi$.
\par
Without the MCF, it is possible to precompute once the solution of the ionic model \eqref{eq:ionic_variables_space_mcf}, appropriately shifting the solution in time by $\psi$ in order to evaluate $w_\mathrm{Ca}$ in \eqref{eq:sarcomere_variables_space_mcf} \cite{stella_fast_2022}. However, when including the MCF, due to its spatially heterogeneous and time dependent nature in \eqref{eq:ionic_variables_space_mcf}, it is not possible to precompute the solutions to \eqref{eq:ionic_variables_space_mcf}, nor is it possible to use the same solution in all points of the domain. Nevertheless, the spatial decoupling of the depolarization mechanism in \eqref{eq:ionic_cable_mcf} allows to significantly reduce computational cost, since the update of the transmembrane potential can be done nodally through \eqref{eq:ionic_cable_mcf}, rather than by solving the system associated to \eqref{eq:monodomain}. Moreover, Equation \eqref{eq:eikonal} allows for the use of much coarser meshes with respect to the monodomain equation \eqref{eq:monodomain}.
\subsection{Numerical approximation}
\label{sec:numerical_approximation}
\begin{figure}[t]
    \centering
    \begin{tikzpicture}[>=latex']
    \tikzset{block/.style= {draw, rectangle, align=center,minimum width=2.8cm,minimum height=1cm},
    }
    \tikzset{blocklarge/.style= {draw, rectangle, align=center,minimum width=7cm,minimum height=2cm},
    }
    \tikzset{blockmiddle/.style= {draw, rectangle, align=center,minimum width=3.5cm,minimum height=1cm},
    }

    \node [blocklarge] (start){};
    \node [block, right = -6.5cm of start]  (A0) {\small
        ionic\\ $\boldsymbol{w}^{n-1,1},u^{n-1,1}$};

    \node [block, right = -3.3cm of start] (A1){\small
        ionic\\ $\boldsymbol{w}^{n},u^{n}$};

    \node [blocklarge, right = 1cm of start] (B2){};
    \node [block, right = -6.5cm of B2] (C1){\small
        ionic\\ $\boldsymbol{w}^{n,1},u^{n,1}$};
    \node [block, right = -3.3cm of B2] (C2){\small
        ionic\\ $\boldsymbol{w}^{n+1},u^{n+1}$};

    \node [blockmiddle, right = -4.5cm of start, above = -4.5cm of start] (D1){\small
        force generation\\ $\boldsymbol{y}^{n}$};
    \node [blockmiddle, right = -4.5cm of B2, above = -4.5cm of B2] (D2){\small
        force generation\\ $\boldsymbol{y}^{n+1}$};

    \node [blockmiddle, below = 1cm of D1] (E1){\small
        mechanics\\ $\mathbf{d}^{n}$};
    \node [blockmiddle, below = 1cm of D2] (E2){\small
        mechanics\\ $\mathbf{d}^{n+1}$};

    \node [blockmiddle, below = 1cm of E1] (F1){\small
        circulation\\ $p^{n}_\mathrm{LV}$};
    \node [blockmiddle, below = 1cm of E2] (F2){\small
        circulation\\ $p^{n+1}_\mathrm{LV}$};

    \node[anchor=south] at (0,1){$t^n$};
    \node[anchor=south] at (8,1){$t^{n+1}$};

    \node [coordinate, below = 1cm of A1] (A){};
    \node [coordinate, left = 1.62cm of A] (B){};
    \node [coordinate, above = 0.3cm of D1.0] (C){};
    \node [coordinate, right = 2.25cm of C] (D){};
    \node [coordinate, above = 2cm of D] (E){};
    \node [coordinate, right = 0.5cm of E] (F){};
    \node [coordinate, below = 1cm of C2.270] (G){};
    \node [coordinate, left = 1.63cm of G] (H){};
    \node [coordinate, above = 0.2cm of E1.0] (I){};
    \node [coordinate, above = 0.2cm of F1.0] (J){};
    \node [coordinate, right = 2.25cm of I] (K){};
    \node [coordinate, right = 2.25cm of J] (L){};
    \node [coordinate, above = 1.62cm of K] (M){};
    \node [coordinate, above = 1.62cm of L] (N){};
    \node [coordinate, right = 2.25cm of M] (O){};
    \node [coordinate, right = 2.25cm of N] (P){};
    \node [coordinate, right = 1.5cm of M] (Q){};
    \node [coordinate, left = 2.5cm of B2.270] (R){};

    \node [coordinate, above = 0.2cm of F1.0] (S){};
    \node [coordinate, right = 2.25cm of S] (T){};
    \node [coordinate, above = 1.62cm of T] (U){};
    \node [coordinate, right = 2.25cm of U] (V){};

    \path[draw,dashed, ->]
        (A0.0) -- (A1.180);
    \path[draw,dashed, ->]
        (C1.0) -- (C2.180);
    \path[draw,dashed, ->]
        (D1.0) -- (D2.180);
    \path[draw,dashed, ->]
        (A1.0) -- (C1.180);
    \path[draw, -]
        (A1.270) -- (A);
    \path[draw, -]
        (B) -- (A) node[midway,below] {\small$w^n_\mathrm{Ca}$};
    \path[draw, ->]
        (B) -- (D1.90);
    \path[draw, -, thick]
        (C) -- (D) ;
    \path[draw, -, thick]
        (D) -- (E) node[midway,right] {\small$y^{n}_\mathrm{dCa_{Tn}}$};
    \path[draw, ->, thick]
        (E) -- (F);
    \path[draw, -]
        (C2.270) -- (G);
    \path[draw, -]
        (G) -- (H)node[midway,below] {\small$w^{n+1}_\mathrm{Ca}$};
    \path[draw, ->]
        (H) -- (D2.90);
    \path[draw,dashed, ->]
        (E1.0) -- (E2.180);
    \path[draw,dashed, ->]
        (F1.0) -- (F2.180);
    \path[draw, ->]
        (D1.270) -- (E1.90) node[midway,right] {\small$T^{n}_\mathrm{a}$};
    \path[draw, ->]
        (E1.270) -- (F1.90) node[midway,right] {\small$V^{n}_\mathrm{LV}$} ;
    \path[draw, ->]
        (D2.270) -- (E2.90) node[midway,right] {\small$T^{n+1}_\mathrm{a}$};
    \path[draw, ->]
        (E2.270) -- (F2.90) node[midway,right] {\small$V^{n+1}_\mathrm{LV}$};
    \path[draw, -]
        (I) -- (K) ;
    \path[draw, -]
        (K) -- (M) node[midway,right] {\small$\mathbf{F}^n$};
    \path[draw, ->]
        (M) -- (O) ;
    \path[draw, ->]
        (Q) -- (R) ;
    \path[draw, -]
        (S) -- (T) ;
    \path[draw, -]
        (T) -- (U) node[midway,right] {\small$p^n_\mathrm{LV}$};
    \path[draw, ->]
        (U) -- (V) ;

    % Calculate the center of the image
    \coordinate (center) at (current bounding box.center);

    % Legend
    \matrix [draw,scale=0.5,anchor=south] at ([yshift=-7.2cm]center) {
        \draw[->] (0,0) -- (0.7,0); & \node[anchor=west,font=\footnotesize] at (1.2,0) {inter-model coupling}; \\[0.07cm]
        \draw[->,thick] (0,0) -- (0.7,0); & \node[anchor=west,font=\footnotesize] at (1.2,0) {mechano-calcium feedback}; \\[0.07cm]
        \draw[->,dashed] (0,0) -- (0.7,0); & \node[anchor=west,font=\footnotesize] at (1.2,0) {intra-model time advancement}; \\
    };
\end{tikzpicture}
    \caption{Representation of the time advancement and coupling scheme of the fully coupled multiscale problem, as in \cite{fedele_comprehensive_2023}. Dashed arrows indicate the time advancement of a single model, while full arrows indicate coupling between models. The thick arrow represents the new coupling \eqref{eq:catn_dt_active_force}.}
    \label{fig:diagram_time_discr}
\end{figure}
In this section we present the numerical scheme suitable for the numerical approximation of the coupling condition \eqref{eq:catn_dt_active_force}. For the numerical approximation of \eqref{eq:full_em_cell}, \eqref{eq:full_em_space} or alternatively \eqref{eq:eik_em_space}, the numerical schemes hinge on what was already presented in \cite{fedele_comprehensive_2023,regazzoni_cardiac_2022}, with the only addition being the numerical treatment of the term \eqref{eq:catn_dt_active_force} in Equations \eqref{eq:ionic_variables_mod}, \eqref{eq:ionic_variables_space} and \eqref{eq:ionic_variables_space_mcf}.
\par
For Equation \eqref{eq:full_em_cell} we use a staggered approach to implement the coupling \eqref{eq:catn_dt_active_force}, with a coarser timescale employed for the force generation model and a finer timescale for the ionic model. Let $\Delta t$ be the time step of the coarse time scale, with $t^n=n\Delta t, \quad N\Delta t = T$, and let it be subdivided further into $N_\mathrm{sub}$ subintervals, such that $\Delta t = N_\mathrm{sub}\tau$, where $\tau$ is the time step of the fine timescale. Let the numerical solutions be denoted as $\boldsymbol{y}^{n}\approx \boldsymbol{y}(t^n)$ for the force generation model and $\boldsymbol{w}^{n,k}\approx \boldsymbol{w}(t^n + k\tau)$ for the ionic model, where in particular it hods that $\boldsymbol{w}^{n,N_\mathrm{sub}}=\boldsymbol{w}^{n+1,0}=\boldsymbol{w}^{n+1}$. Then, by employing a backward Euler scheme of order one for the approximation of the time derivative, the solution scheme for \eqref{eq:ionic_variables_mod} is the following implicit-explicit scheme \cite{regazzoni_cardiac_2022}, for all $k=0,\dots,N_\mathrm{sub}-1,\ n=0,\dots,N-1$:
\begin{equation}
\label{eq:ionic_solution_scheme}
    \frac{\boldsymbol{w}^{n,k+1}-\boldsymbol{w}^{n,k}}{\tau} = \tilde{\boldsymbol{g}}\left(u^{n,k}, \boldsymbol{w}^{n,k+1}_\mathrm{G}, \boldsymbol{w}^{n,k}_\mathrm{I}, y_{\mathrm{dCa_{Tn}}}^n\right), 
\end{equation}
which is implicit for the gating variables $\boldsymbol{w}_\mathrm{G}$, explicit for the ionic concentrations $\boldsymbol{w}_\mathrm{I}$. In particular, scheme \eqref{eq:ionic_solution_scheme} has the following form:
\begin{subequations}
    \label{eq:imex_extended}
    \begin{empheq}[left=\empheqlbrace]{align}
    & \frac{\boldsymbol{w}^{n,k+1}_\mathrm{G}-\boldsymbol{w}^{n,k}_\mathrm{G}}{\tau} = \tilde{\boldsymbol{g}}_\mathrm{G}\left(u^{n,k}, \boldsymbol{w}^{n,k+1}_\mathrm{G}\right),\\
    & \frac{\boldsymbol{w}^{n,k+1}_\mathrm{I}-\boldsymbol{w}^{n,k}_\mathrm{I}}{\tau} = \tilde{\boldsymbol{g}}_\mathrm{I}\left(u^{n,k}, \boldsymbol{w}^{n,k+1}_\mathrm{G}, \boldsymbol{w}^{n,k}_\mathrm{I}, y_{\mathrm{dCa_{Tn}}}^n\right).
    \end{empheq}
\end{subequations}
The numerical approximation of the feedback term corresponds to using a first-order extrapolation on the coarse scale of the calcium-troponin buffering rate \eqref{eq:derived_ca_tn_active_force}. The numerical treatment of the coupling is analogous in the case of the multiscale problems \eqref{eq:full_em_space} and \eqref{eq:eik_em_space}, and is briefly summarized in Figure \ref{fig:diagram_time_discr}, although special care must be employed in the definition of space-dependent quantities, as we will see.
\par
The multiscale problems \eqref{eq:full_em_space} and \eqref{eq:eik_em_space} are solved using continuous nodal finite elements (FE) \cite{fedele_comprehensive_2023,quarteroni_numerical_2017}. In particular, problem \eqref{eq:full_em_space} is solved on two meshes, a coarse mesh $\mathcal{T}_H$ with the FE spaces ${V_H}^{\boldsymbol{y}}$ and ${V_H}^{\mathbf{d}}$ (of equal degree) for Equations \eqref{eq:sarcomere_variables_space} and \eqref{eq:mechanics}, and a fine mesh $\mathcal{T}_h$ with the FE spaces ${V_h}^{\boldsymbol{w}}$ and ${V_h}$ (of equal degree) for Equations \eqref{eq:ionic_variables_space} and \eqref{eq:monodomain}, with the superscript indicating the approximated variables. Systems \eqref{eq:ionic_variables_space} and \eqref{eq:sarcomere_variables_space} are solved nodally \cite{krishnamoorthi_numerical_2013}. In this case, denoting with the subscripts $h$ and $H$ the numerical solutions in their respective approximation spaces, and remembering that $\boldsymbol{w}^{n,k}\in {V_h}^{\boldsymbol{w}}$ and $y_{\mathrm{dCa_{Tn}}}^n\in {V_H}$, the scheme \eqref{eq:ionic_solution_scheme} becomes:
\begin{equation}
\label{eq:feedback_intergrid_ms}
    \frac{\boldsymbol{w}^{n,k+1}_h-\boldsymbol{w}^{n,k}_h}{\tau} = \tilde{\boldsymbol{g}}\left(u^{n,k}_h, \boldsymbol{w}^{n,k+1}_{\mathrm{G},h}, \boldsymbol{w}^{n,k}_{\mathrm{I},h}, y_{\mathrm{dCa_{Tn}},H,h}^n\right),
\end{equation}
where $y_{\mathrm{dCa_{Tn}},H,h}^n$ is the interpolation of $y_{\mathrm{dCa_{Tn}},H}^n\in V_H$ in $V_h$ obtained via radial basis functions \cite{salvador_intergrid_2020,bucelli_robust_2024}. The rest of the system is solved as previously presented in \cite{fedele_comprehensive_2023,regazzoni_cardiac_2022,bucelli_preserving_2023}, by exploiting intergrid-staggered algorithms and radial basis function interpolation.
\par
As for Equation \eqref{eq:eik_em_space}, due to less restrictive mesh resolution requirements \cite{colli_franzone_wavefront_1990}, the problem is solved on a single coarse level mesh $\mathcal{T}_h$ and the same polynomial degree is used for all approximation spaces. In this case the feedback condition in space is trivially obtained by replacing $y_{\mathrm{dCa_{Tn}},H,h}^n$ with $y_{\mathrm{dCa_{Tn}},h}^n$ directly in \eqref{eq:feedback_intergrid_ms}. Finally, we remark once again that the solution of the ionic model \eqref{eq:ionic_variables_space_mcf} cannot be done offline and once-for-all as in \cite{stella_fast_2022}, due to the time dependent and spatially heterogeneous nature of the feedback condition.
\section{Numerical results}
\label{sec:numerical_results}
In this section we present numerical results for the models introduced in the previous section. The approach presented in Section \ref{sec_model_derivation} generally applies to a wide range of families of ionic and sarcomere models. The principal aim of this section is understanding the implications of using different formulations of the previously introduced electromechanical models in Sections \ref{sec:fully_coupled_multiscale} and \ref{sec:eikonal_driven} by comparing them starting from similar baseline cellular electromechanical models.
\par
In this work, we consider the ten Tusscher and Panfilov (TTP06) \cite{ten_tusscher_alternans_2006} ionic model due to its widespread use, coupled with the RDQ20 model as the basis for our cellular electromechanical coupling. In Section \ref{sec:cellular_electromechanical_model} we describe how the approach introduced in this paper is applied in this case, along with a mathematically and biophysically consistent calibration procedure to fit its new parameters. In Section \ref{sec:numerical_results_multiscale} we present and compare numerical results for the multiscale models presented in Sections \ref{sec:fully_coupled_multiscale} and \ref{sec:eikonal_driven} starting from the  cellular model presented in Section \ref{sec:cellular_electromechanical_model}. In Section \ref{sec:numerical_calcium_sens} we examine the effects of parameter perturbation on the models presented in Sections \ref{sec:cellular_electromechanical_model} and \ref{sec:numerical_results_multiscale} to gain further insight into the multiscale effects of the different formulations.
\par
Models and numerical schemes were implemented in the software library \texttt{life\textsuperscript{x}} \cite{africa_life_2022}, which is based on the finite-element framework offered by \texttt{dealii} \cite{arndt_dealii_2021,africa_dealii_2024}. Multiscale simulations were ran on 72 cores on the LEONARDO supercomputer (two Intel Sapphire Rapids CPUs at 2.00 GHz, 512 GB RAM) in the CINECA supercomputing center (Italy). Single-cell simulations were ran on a single core on a personal desktop computer (Intel Core i5-9600K CPU). All numerical and model parameters are reported in Appendix \ref{app:model_numerical_parameters}.

\subsection{Cellular electromechanical model calibration}
\label{sec:cellular_electromechanical_model}
We recall that the TTP06 model is of the single calcium buffer type, thus following Section \ref{sec:single_buffer} we apply split \eqref{eq:split_buffer}, introducing a free parameter ${\mathrm{Tn}}_{\mathrm{c,max}}$ in the ionic model. The introduction of the feedback mechanism affects both the intracellular calcium transient and the active force transient of the cellular electromechanical model, both depending on the free parameter, but also other parameters of the ionic and active force models. This leads, as discussed in Section \ref{sec:single_buffer}, to a need to recalibrate the model \eqref{eq:full_em_cell} not limited to the free parameter ${\mathrm{Tn}}_{\mathrm{c,max}}$.
\par
In order to recover consistent calcium transients and contraction kinetics while incorporating the feedback, we tune the new model, henceforth denoted with TTP06+RDQ20, and perform a two-step calibration of the cellular model \eqref{eq:full_em_cell}, where the first step is concerned with recovering the initial calcium transients \cite{ten_tusscher_alternans_2006} through the tuning of the ionic parameters, and the second part concerns the recovery of realistic contraction kinetics \cite{regazzoni_biophysically_2020, fedele_comprehensive_2023}. Unlike previous studies \cite{margara_-silico_2021, mazhar_detailed_2024}, we adjust the electrophysiological parameters to align with the original model and not necessarily just the experimental data \cite{coppini_late_2013}. In this way we ensure that we maintain consistent baseline cellular electromechanical behavior across the cellular models, and avoid introducing discrepancies in single-cell behaviors that could mask multiscale-level effects of the mechano-calcium feedback.
\par
{We recast the calibration problem as finding the parameters $\left(\theta_\mathrm{ion},\boldsymbol{\theta}_\mathrm{af}\right)\in R_1\times \boldsymbol{R}_2$ as solutions to a particular minimization problem. The parameters $\theta_\mathrm{ion}$ are the ionic parameters tuned during the first step, which we will consider to be a single parameter, while $\boldsymbol{\theta}_\mathrm{af}=\left(\boldsymbol{\theta}_\mathrm{af,kin},a_\mathrm{XB}\right)$ are the tunable parameters of the force generation model and are split into the kinetic parameters $\boldsymbol{\theta}_\mathrm{af,kin}=\left(k_\mathrm{basic}, k_\mathrm{off}\right)$ and the upscaling constant $a_\mathrm{XB}$. The admissible ranges for $\theta_\mathrm{ion}$ and $\boldsymbol{\theta}_\mathrm{af}$ are given as $R_1\subset\subset\mathbb{R}$ and $\boldsymbol{R}_2=\boldsymbol{R}_{2,\mathrm{kin}}\cross\mathbb{R}^+$, with $\boldsymbol{R}_{2,\mathrm{kin}}\subset\subset\mathbb{R}^2$.}
\par
{For the solution of the minimization problem, we will consider numerical approximations of the limit cycle solutions of system \eqref{eq:full_em_cell} in normal physiological conditions (e.g. without alternans) and at constant sarcomere lengths, obtained by applying periodic stimuli of type \eqref{eq:applied_current}. The local minima are computed using a
grid search algorithm. Due to the time-discrete nature of the approximated solutions, any continuous limit cycle $x(t)$ approximated using a timestep equal to $\Delta t$ reduces to a finite-dimensional vector $\boldsymbol{x}$ such that:}
\begin{equation*}
    \boldsymbol{x}=\left(x^0,x^1,\dots, x^{N}\right),\quad x^k\approx x(k\Delta t), \quad k = 0,\dots, N.
\end{equation*}
We write the general calibration problem as follows.
\par
\textit{{Let $\Delta t>0$ be a fixed time discretization step, and let $\boldsymbol{w}^0_\mathrm{Ca}$ and $\boldsymbol{T}^0_\mathrm{a}$ be the corresponding limit cycle approximations of \eqref{eq:full_em_cell} without the feedback mechanism \eqref{eq:catn_dt_active_force}, corresponding to a vector with $N+1$ components, for a set of fixed ionic and activation model parameters $\hat{\theta}_\mathrm{ion},\hat{\boldsymbol{\theta}}_\mathrm{af,kin},\hat{a}_\mathrm{XB}$. Let $\boldsymbol{w}_\mathrm{Ca}(\Tilde{\theta}_\mathrm{ion})$ and $\boldsymbol{T}_\mathrm{a}(\Tilde{\theta}_\mathrm{ion},\Tilde{\boldsymbol{\theta}}_\mathrm{af,kin},\Tilde{a}_\mathrm{XB})$ denote limit cycle approximations of \eqref{eq:full_em_cell} including the feedback mechanism \eqref{eq:catn_dt_active_force} for $\Delta t>0$, obtained for the parameters $\Tilde{\theta}_\mathrm{ion},\Tilde{\boldsymbol{\theta}}_\mathrm{af,kin},\Tilde{a}_\mathrm{XB}$. Moreover, let us denote with $T_\mathrm{a,max}$ and $T^0_\mathrm{a,max}$ the approximated maximal active tensions, i.e.
\begin{equation*}
\begin{aligned}
    T_\mathrm{a,max}(\Tilde{\theta}_\mathrm{ion},\Tilde{\boldsymbol{\theta}}_\mathrm{af,kin},\Tilde{a}_\mathrm{XB}) &= \underset{k=0,\dots,N}{\max}{T}^{k}_\mathrm{a}(\Tilde{\theta}_\mathrm{ion},\Tilde{\boldsymbol{\theta}}_\mathrm{af,kin},\Tilde{a}_\mathrm{XB}),
    \\ T^0_\mathrm{a,max}&=\underset{k=0,\dots,N}{\max}{{T}^{0,k}_\mathrm{a}}.
\end{aligned}
\end{equation*}
\newline
Find $\theta_\mathrm{ion}\in R_1$, ${\boldsymbol{\theta}}_\mathrm{af,kin}\in\boldsymbol{R}_{2,\mathrm{kin}}$, $a_\mathrm{XB}\in\mathbb{R}^+$, such that:
\begin{equation}
\label{eq:theta_ion_def}
    \theta_\mathrm{ion} = \underset{\Tilde{\theta}\in R_1}{\argmin} \varphi_\mathrm{ion}(\Tilde{\theta}),
\end{equation}
\begin{equation}
\label{eq:theta_afkin_def}
\boldsymbol{\theta}_\mathrm{af,kin} = \underset{\Tilde{\boldsymbol{\theta}}_\mathrm{af,kin}\in \boldsymbol{R}_{2,\mathrm{kin}}}{\argmin} \varphi_\mathrm{af,kin}(\Tilde{\boldsymbol{\theta}}_\mathrm{af,kin}),
\end{equation}
and
\begin{equation}
\label{eq:axb_calib}
    a_\mathrm{XB} =  \frac{T^0_\mathrm{a,max}}{T_\mathrm{a,max}(\theta_\mathrm{ion},{\boldsymbol{\theta}}_\mathrm{af,kin},\hat{a}_\mathrm{XB})}\Hat{a}_\mathrm{XB},
\end{equation}
where:
\begin{equation}
\label{eq:cost_ion}
    \varphi_\mathrm{ion}(\Tilde{\theta}) = |\boldsymbol{w}^0_\mathrm{Ca} - \boldsymbol{w}_\mathrm{Ca}(\Tilde{\theta})|,
\end{equation}
and
\begin{equation}
\label{eq:cost_af}
    \varphi_\mathrm{af,kin}(\Tilde{\boldsymbol{\theta}}_\mathrm{af,kin}) = \Bigg|\frac{\boldsymbol{T}^0_\mathrm{a}}{T^0_\mathrm{a,max}} - \frac{\boldsymbol{T}_\mathrm{a}(\theta_\mathrm{ion},\Tilde{\boldsymbol{\theta}}_\mathrm{af,kin},\hat{a}_\mathrm{XB})}{T_\mathrm{a,max}(\theta_\mathrm{ion},\Tilde{\boldsymbol{\theta}}_\mathrm{af,kin},\hat{a}_\mathrm{XB})}\Bigg|,
\end{equation}
where $|\cdot|$ is the Euclidean norm.}
}
\par
{If we assume $\boldsymbol{w}_\mathrm{Ca}$ to have continuous dependence on $\Tilde{\theta}_\mathrm{ion}$ in $R_1$, and $\boldsymbol{T}_\mathrm{a}$ to have continuous dependence on $\Tilde{\boldsymbol{\theta}}_\mathrm{af,kin}$ in $\boldsymbol{R}_{2,\mathrm{kin}}$, then both $\theta_\mathrm{ion}$ and ${\boldsymbol{\theta}}_\mathrm{af,kin}$ defined in \eqref{eq:theta_ion_def} and \eqref{eq:theta_afkin_def} exist, as consequence of the Weierstrass theorem. This calibration problem can be reinterpreted as, given a parametrization of a cellular electromechanical model and a certain discrepancy measure, find a best fit in another family of cellular models according to this discrepancy measure.}
\par
{In particular, we consider the initial cellular electromechanical model to be the TTP06+RDQ20 model without the MCF, the initial parametrization to be one yielding realistic ventricular contraction properties \cite{fedele_comprehensive_2023}, the other family of cellular electromechanical models to be the TTP06+RDQ20 model with the MCF, and the discrepancy being measured in terms of intracellular calcium transients \eqref{eq:cost_ion} and active force kinetics \eqref{eq:cost_af}. From the form of the calibration problem, we see that, in fact, it consists of two sequential minimization problems, the first one recovering the intracellular calcium transients, and the second one recovering active force kinetics. Although the ionic parameter $\theta_\mathrm{ion}$ is recalibrated in order to minimize the changes in the calcium transient curve, due to the shift of the calcium transient waveform, the force generation parameters $\boldsymbol{\theta}_\mathrm{af}$ need to be recalibrated in order to capture the baseline active force waveform. In particular, the active force kinetics is regulated through $\boldsymbol{\theta}_\mathrm{af,kin}=\left(k_\mathrm{basic}, k_\mathrm{off}\right)$, while the active force amplitude is regulated through the upscaling constant $a_\mathrm{XB}$ \cite{regazzoni_biophysically_2020}. For the single-cell model \eqref{eq:full_em_cell} the active force amplitude has no impact on the state variables' dynamics, therefore the calibration of the kinetic parameters can be done separately by minimizing $\varphi_\mathrm{af,kin}$.}
\par
{In principle, the active force $\boldsymbol{T}_\mathrm{a}$ depends also on the upscaling constant $\hat{a}_\mathrm{XB}$ by a multiplicative factor. However by dividing it by $T_\mathrm{a,max}$ in \eqref{eq:cost_af} that dependence is lost, making the function $\varphi_\mathrm{af,kin}$ independent of $\hat{a}_\mathrm{XB}$, and therefore the choice of ${a}_\mathrm{XB}$ after minimizing $\varphi_\mathrm{ion}$ and $\varphi_\mathrm{af,kin}$ still arbitrary. By requiring that the reparametrized model preserves the maximal active force, the new crossbridge stiffness $a_\mathrm{XB}$ can be computed with respect to its original value $\Hat{a}_\mathrm{XB}$ through relation \eqref{eq:axb_calib}.}
\par
{We remark that the choice of the parameter $\theta_\mathrm{ion}$ is not straightforward. Indeed, the choice of the ranges $R_1$ and $\boldsymbol{R}_{2,\mathrm{kin}}$ is not unique, and therefore in the case where the minima of $\varphi_\mathrm{ion}$ and $\varphi_\mathrm{af,kin}$ are located on the boundaries of $R_1$ and $\boldsymbol{R}_{2,\mathrm{kin}}$ the calibration may yield ambiguous and uninterpretable results.}
\par
\begin{figure}[H]
     \centering
     \subfloat[][]{\includegraphics[width = 0.48\textwidth]{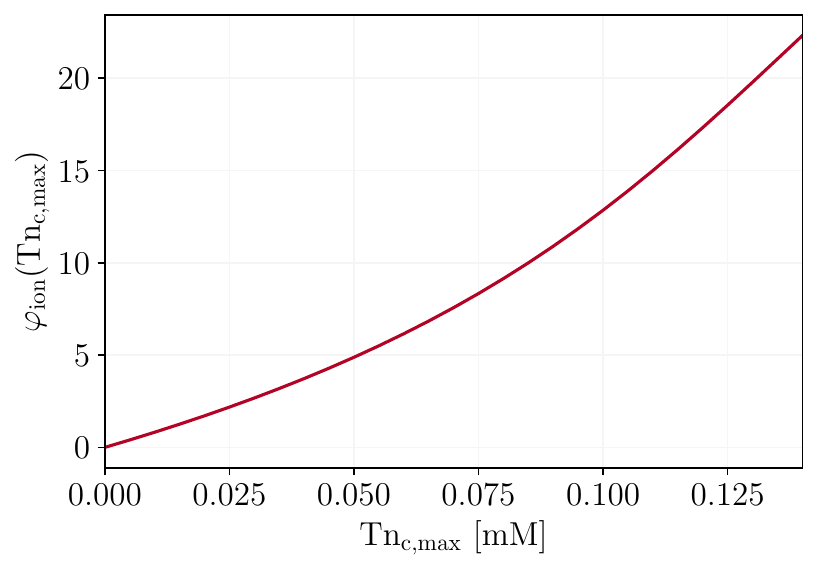}\label{fig:costs_vs_tnc}}\hfill
     \subfloat[][]{\includegraphics[width = 0.48\textwidth]{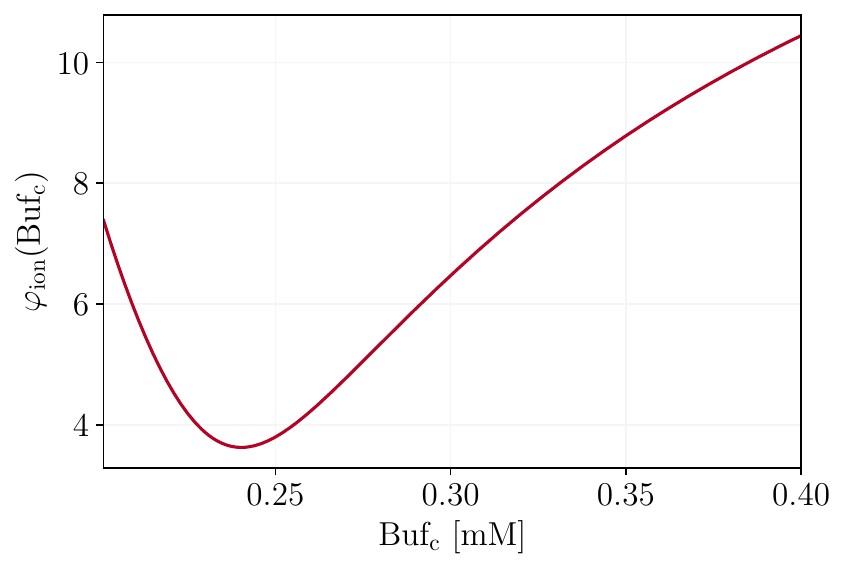}\label{fig:costs_vs_bufc}}
     \caption{Value of the discrepancy metric $\varphi_\mathrm{ion}$ depending on the choice of the parameter $\theta_\mathrm{ion}$. (a): Discrepancy metric $\varphi_\mathrm{ion}$ depending of the maximal calcium-troponin buffer concentration $\mathrm{Tn}_\mathrm{c,max}$. (b): Discrepancy metric $\varphi_\mathrm{ion}$ depending of the total cytosolic buffer concentration $\mathrm{Buf}_\mathrm{c}$.}
\end{figure}

In fact, given the remarks in Section \ref{sec:single_buffer}, that the model including the feedback with ${\mathrm{Tn}}_{\mathrm{c,max}}=0$ is equivalent to the model without the feedback, if the calibration parameter $\theta_\mathrm{ion}$ is taken to be ${\mathrm{Tn}}_{\mathrm{c,max}}$, it would have a degenerate minimum and would not yield other local minima within physiological ranges, as exemplified in Figure \ref{fig:costs_vs_tnc}. Instead, we set ${\mathrm{Tn}}_{\mathrm{c,max}}$ consistently with literature values \cite{smith_calcium_2019}, and calibrate with respect to $\overline{\mathrm{Buf}}_{\mathrm{c,max}}$, obtaining the minimum point depicted in Figure \ref{fig:costs_vs_bufc}.
\par
\begin{figure}[H]
    \centering
    \begin{minipage}[c]{0.48\textwidth}
        \centering
        \includegraphics[width=\textwidth]{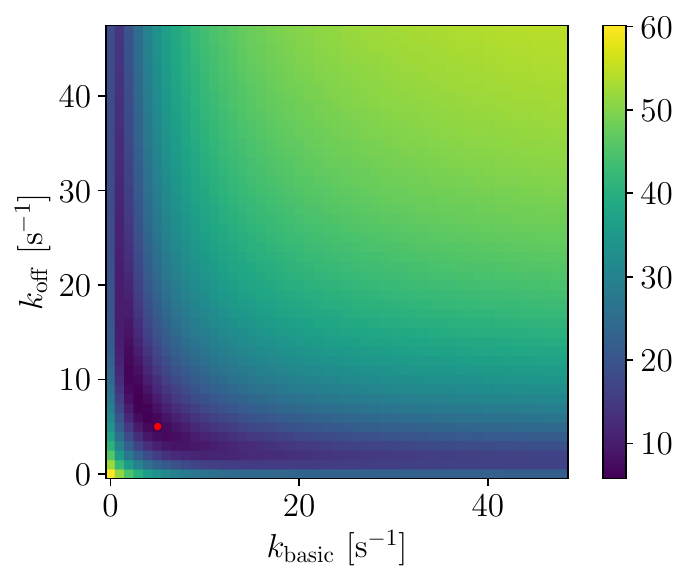}
        \caption{Discrepancy metric $\varphi_\mathrm{af,kin}$ depending of the kinetic parameters $k_\mathrm{off},k_\mathrm{basic}$ in $\boldsymbol{R}_{2,\mathrm{kin}}$. The red dot represents the parameters at which the minimum is attained.}
        \label{fig:kinetic_costs_heatmap}
    \end{minipage}\hfill
    \begin{minipage}[c]{0.48\textwidth}
        \centering
        \vspace{4mm}
        \includegraphics[width=\textwidth]{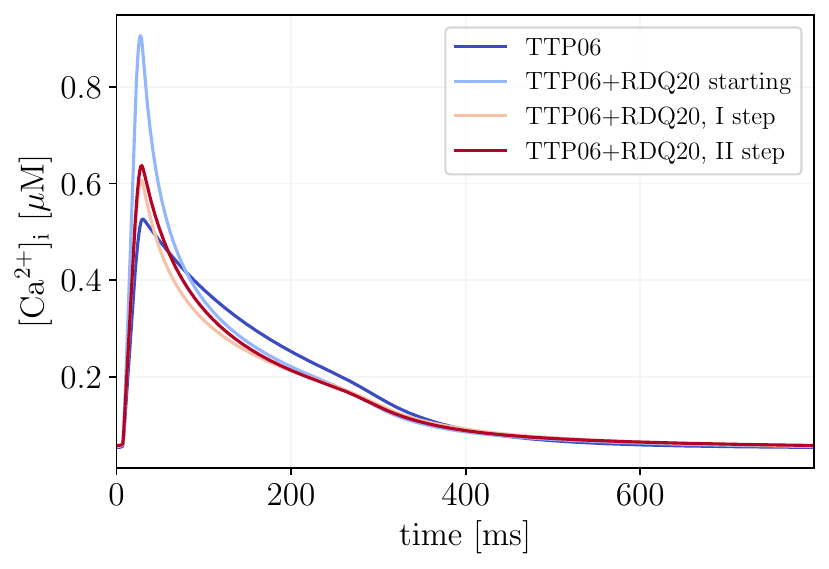}
        \vspace{0mm}
        \caption{Free intracellular calcium $[\mathrm{Ca}^{2+}]_\mathrm{i}$ for each step of the calibration process.}
        \label{fig:calcium_calibration}
    \end{minipage}
\end{figure}

With $\theta_\mathrm{ion}$ fixed in the previous step, the minimization of \eqref{eq:cost_af} as well yields a local minimum in the interior of the parameter space $\boldsymbol{R}_{2,\mathrm{kin}}$, depicted in Figure \ref{fig:kinetic_costs_heatmap}. Finally with $\theta_\mathrm{ion}$ and $\boldsymbol{\theta}_\mathrm{af,kin}$ determined by the two previous steps, we fix the crossbridge stiffness as in \eqref{eq:axb_calib}.
\par
\begin{figure}[t]
     \centering
     \subfloat[][]{\includegraphics[width = 0.48\textwidth]{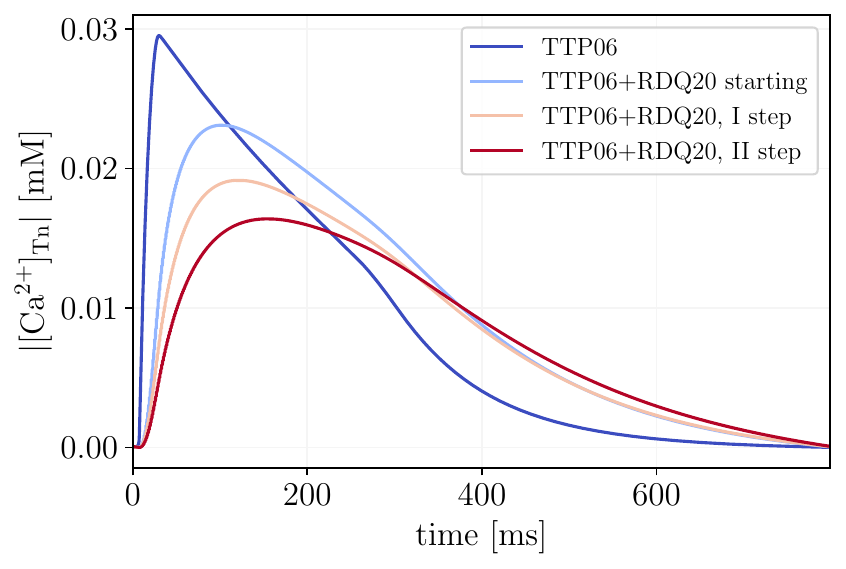}\label{fig:catrpn}}\hfill
     \subfloat[][]{\includegraphics[width = 0.48\textwidth]{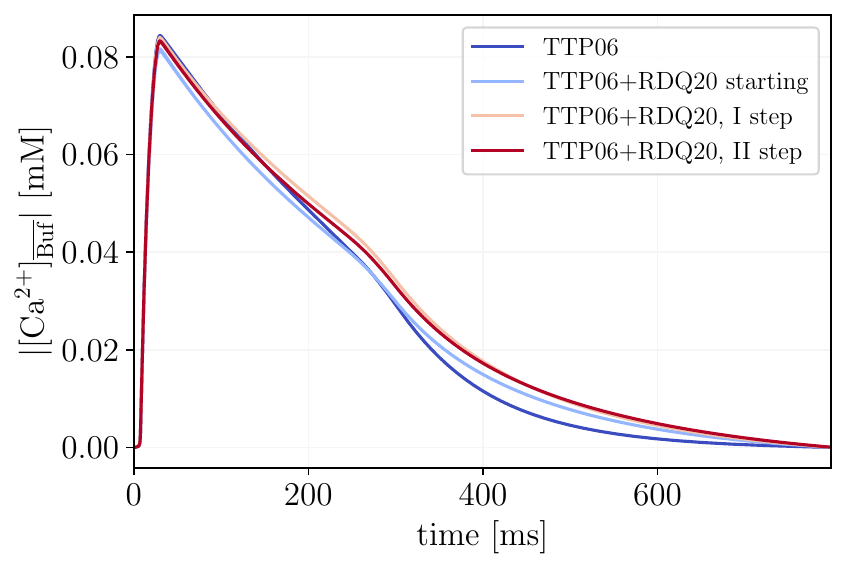}\label{fig:total_ca_buf}}
     \caption{Buffer-bound calcium concentration amplitudes for each step of the calibration process. (a): Troponin-bound calcium concentration $[\mathrm{Ca}^{2+}]_\mathrm{Tn}$ amplitude. (b): Total buffer-bound (troponin and other) calcium concentration $[\mathrm{Ca}^{2+}]_{\overline{\mathrm{Buf}}}$ amplitude.}
     \label{fig:buffers}
\end{figure}

\begin{figure}[t]
     \centering
     \subfloat[][]{\includegraphics[width = 0.48\textwidth]{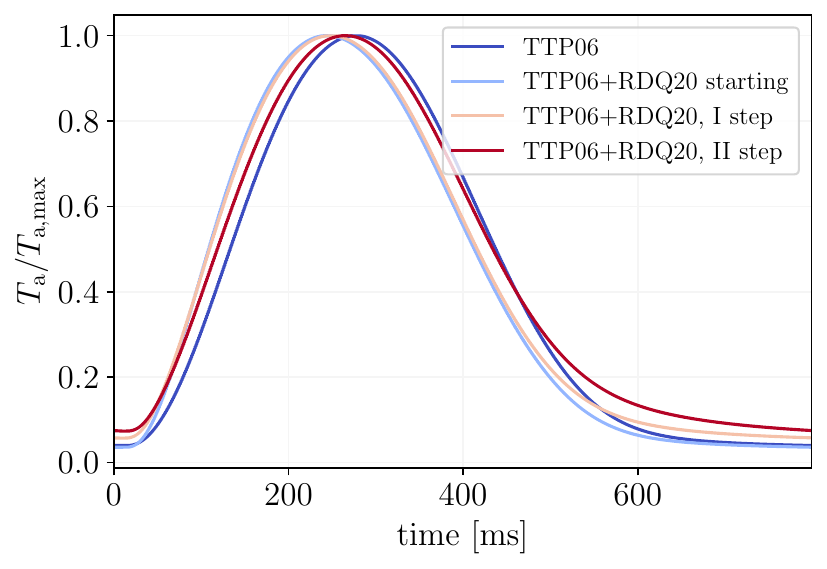}\label{fig:af_kinetics}}\hfill
     \subfloat[][]{\includegraphics[width = 0.48\textwidth]{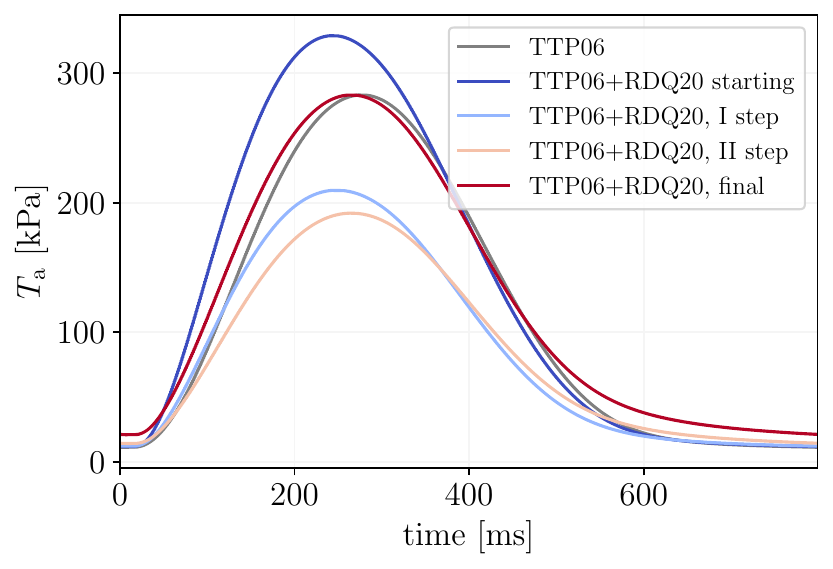}\label{fig:af_calibration}}
     \caption{Active force transients for each step of the calibration process. (a): Normalized active force $T_\mathrm{a}/T_\mathrm{a,max}$. (b): Active force $T_\mathrm{a}$.}
     \label{fig:active_force_calibrations}
\end{figure}

The results following each step of calibration are reported in Figures \ref{fig:calcium_calibration}, \ref{fig:buffers} and \ref{fig:active_force_calibrations}. In Figure \ref{fig:calcium_calibration} we can see that the initial, uncalibrated TTP06+RDQ20 model has a significantly higher peak calcium concentration with respect to the TTP06 model. This is due to the fact that the inclusion of the MCF condition slows down the buffering of the free cytosolic calcium $[\mathrm{Ca}^{2+}]_\mathrm{i}$ following its sudden influx, owing to the small reaction rate constants $k_\mathrm{off}, k_\mathrm{basic}$ associated to the calcium-troponin buffering rate, as demonstrated by the modification in the troponin-bound calcium waveform reported in Figure \ref{fig:catrpn} and discussed in Section \ref{sec:feedback_condition}. The calibration of $\varphi_\mathrm{ion}$ with respect to $\overline{\mathrm{Buf}}_{\mathrm{c,max}}$ can be interpreted as solving the problem of an excessive free calcium peak by adding more of the fast generic buffer $\mathrm{Buf}$, whose effect can be seen on the total buffer-bound concentration in Figure \ref{fig:total_ca_buf}. This beneficial effect is counteracted by the lowering of the free cytosolic calcium during the relaxation stage of the calcium transient waveform owing to the higher free calcium capture by the combined buffers, seen in the relaxation part of the total buffer-bound calcium wave in Figure \ref{fig:total_ca_buf}.

Lastly, the effects of the calibration on the active force transients are shown in Figure \ref{fig:af_calibration}, and although the shape of the active force transient remains roughly similar, the newly calibrated bidirectionally coupled model exhibits a higher value of the diastolic active force, possibly owed to higher diastolic troponin-bound calcium concentrations. The changes to the characteristic times of the intracellular calcium and active force transient waves are reported in Table \ref{table:cat_biomarkers_computed_1}, where we can see that the inclusion of the MCF delays slightly the overall relaxation time of the active force transients, which may translate into negative relaxation properties in the multiscale model. For the free intracellular calcium, an opposite effect is achieved, where although the $50\%$ relaxation time is prolonged, possibly owing to the overall higher amplitude to be relaxed, the $90\%$ relaxation time is reduced.
\begin{table}[H]
\centering
    \begin{tabular}{ p{9em} | c | c }
    Characteristic time & TTP06+RDQ20 no MCF &  TTP06+RDQ20 + MCF\\ 
    \hline\hline
     $t_\mathrm{tp}([\mathrm{Ca}^{2+}]_\mathrm{i})$ &  30.10 $\si{\milli \second}$ & 28.70 $\si{\milli \second}$\\
     $\mathrm{rt}_{50}([\mathrm{Ca}^{2+}]_\mathrm{i})$  & 30.15 $\si{\milli \second}$ & 40.45 $\si{\milli \second}$ \\ 
     $\mathrm{rt}_{90}([\mathrm{Ca}^{2+}]_\mathrm{i})$ & 303.45 $\si{\milli \second}$ & 277.45 $\si{\milli \second}$ \\
     $t_\mathrm{tp}(T_\mathrm{a})$ &  276.00 $\si{\milli \second}$ & 265.00 $\si{\milli \second}$ \\
     $\mathrm{rt}_{50}(T_\mathrm{a})$ & 764.00 $\si{\milli \second}$ & 785.00 $\si{\milli \second}$ \\ 
     $\mathrm{rt}_{90}(T_\mathrm{a})$ & 791.00 $\si{\milli \second}$ & 796.00 $\si{\milli \second}$ \\
    \end{tabular}
    \caption{Time to peak ($t_\mathrm{tp}$), 50\% relaxation time ($\mathrm{rt}_{50}$), and 90\% relaxation time ($\mathrm{rt}_{90}$) for the calcium and active force transient waveforms, for the calibrated TTP06+RDQ20 model with the MCF (TTP06+RDQ20 + MCF), and for the original TTP06 and RDQ20 models without the feedback (TTP06+RDQ20 no MCF).}
    \label{table:cat_biomarkers_computed_1}
\end{table}
\subsection{Results of the multiscale model}
\label{sec:numerical_results_multiscale}
{In this section we study the impact of the inclusion of the MCF in multiscale models, by comparing the introduced cardiac electromechanics models, M+MCF and E+MCF, against each other and their counterparts obtained by neglecting the feedback \eqref{eq:derived_ca_tn_active_force}, simply denoted by M and E.} In order to have similar cellular electromechanical behaviors, for the multiscale models M+MCF and E+MCF we use the TTP06+RDQ20 model for the cellular electromechanics with the calibration obtained at the end of Section \ref{sec:cellular_electromechanical_model}, otherwise, for the M and E models, we use the uncoupled TTP06 and RDQ20 models with the original parameter calibration. Prior to the start of each simulation, the single-cell models were ran for 1000 cycles in order to have initial conditions close to a limit cycle.
\par
\begin{figure}[H]
     \centering
     \subfloat[][]{\includegraphics[width = 0.48\textwidth]{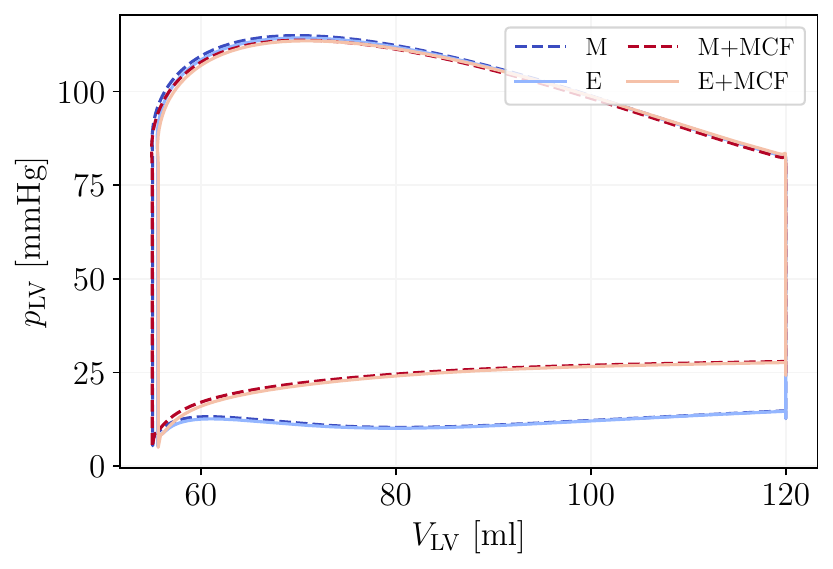}\label{fig:pv_multiscale}}\hfill
     \subfloat[][]{\includegraphics[width = 0.48\textwidth]{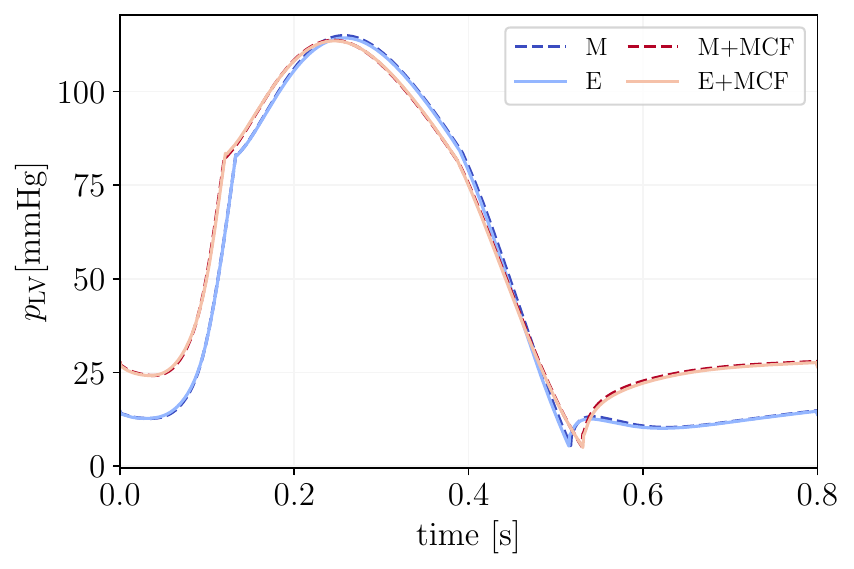}\label{fig:p_lv}}
     \caption{PV loop and endocardial pressure obtained from the multiscale models, both for the full and eikonal-driven problems, both with and without the MCF. (a): PV loop. (b): Endocardial pressure $p_\mathrm{LV}$ trace.}
     \label{fig:results_multiscale_aggregate}
\end{figure}
\begin{figure}[H]
     \centering
     \subfloat[][]{\includegraphics[width = 0.32\textwidth]{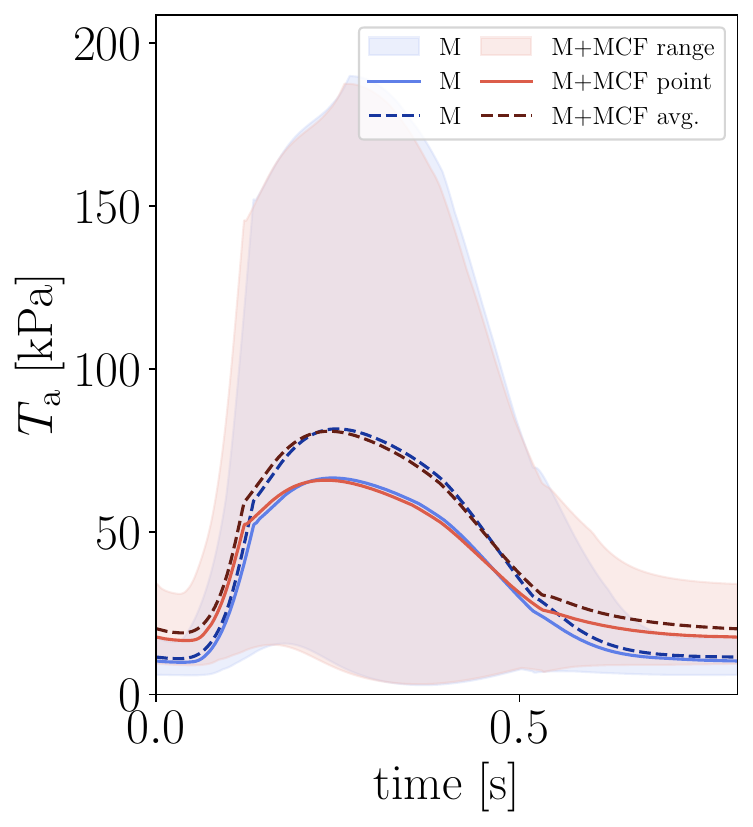}\label{fig:as_ss}}\hfill
     \subfloat[][]{\includegraphics[width = 0.32\textwidth]{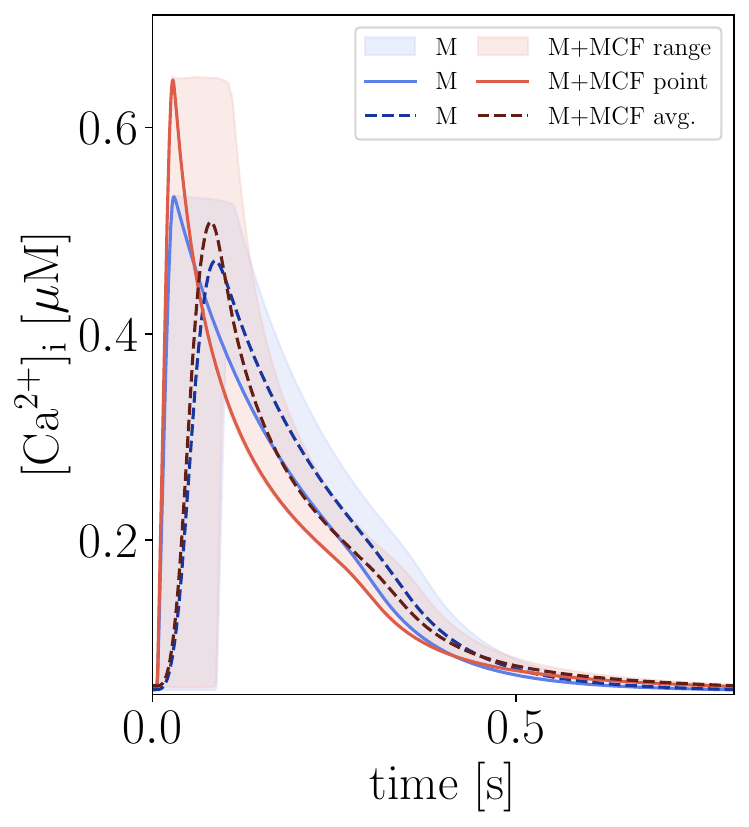}\label{fig:intraca_ss}}\hfill
     \subfloat[][]{\includegraphics[width = 0.32\textwidth]{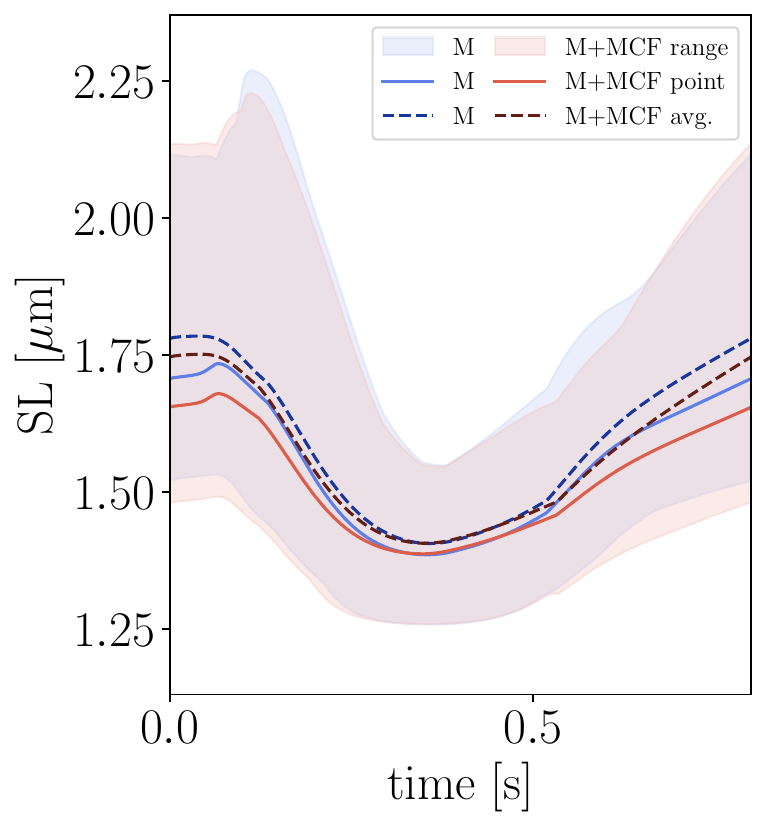}\label{fig:sl_ss}}
     \caption{Space-dependent quantities comparison between the M and M+MCF models (range of values, pointwise value, average value) for the active stress, free intracellular calcium concentration and sarcomere lengths. (a): Active stress $T_\mathrm{a}$ space statistics. (b): Free intracellular calcium concentration $[\mathrm{Ca}^{2+}]_\mathrm{i}$ space statistics. (c): Sarcomere length $\mathrm{SL}$ space statistics.}
     \label{fig:spatial_statistics}
\end{figure}
\begin{figure}[H]
     \centering
     \subfloat[][]{\includegraphics[width = 0.32\textwidth]{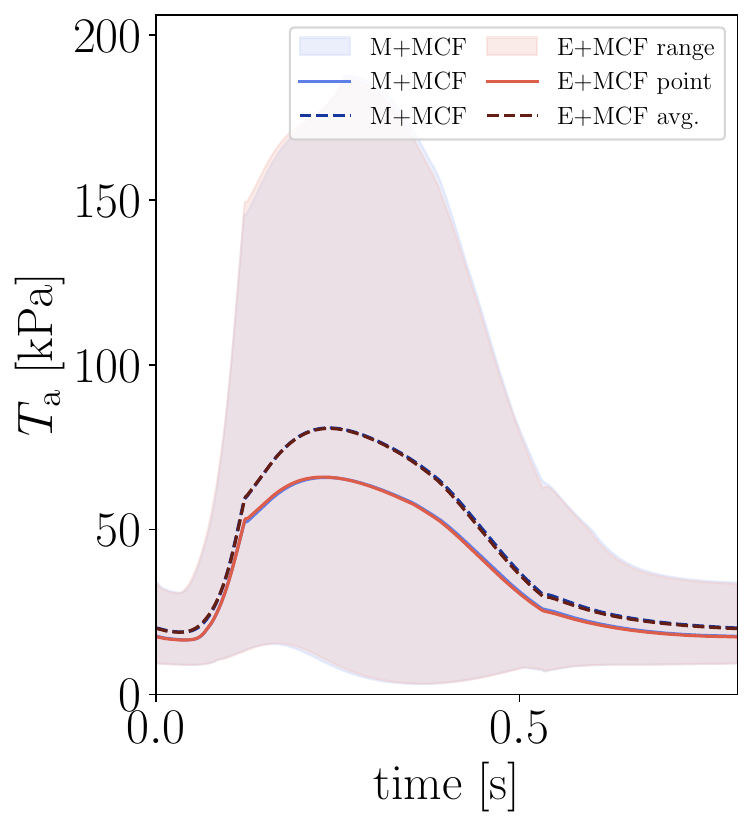}\label{fig:as_ss_mmem}}\hfill
     \subfloat[][]{\includegraphics[width = 0.32\textwidth]{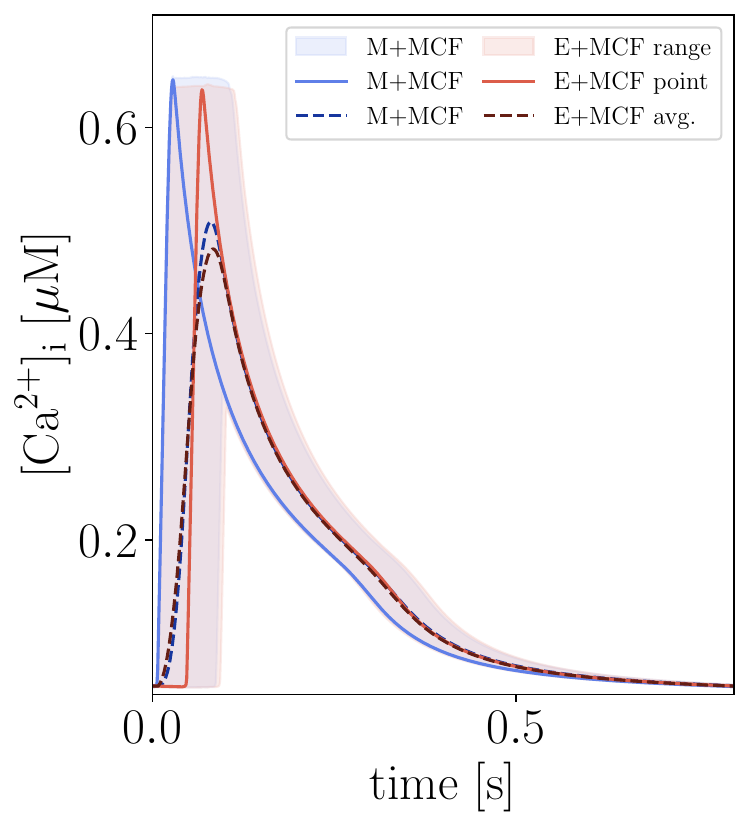}\label{fig:intraca_ss_mmem}}\hfill
     \subfloat[][]{\includegraphics[width = 0.32\textwidth]{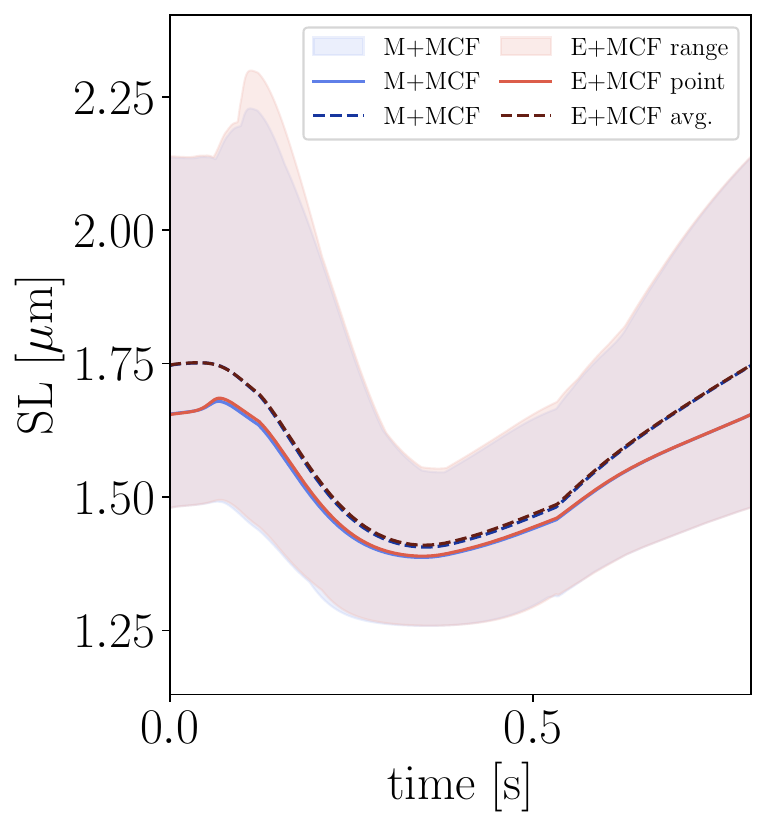}\label{fig:sl_ss_mmem}}
     \caption{Space-dependent quantities comparison between the M+MCF and E+MCF models (range of values, pointwise value, average value) for the active stress, free intracellular calcium concentration and sarcomere lengths. (a): Active stress $T_\mathrm{a}$ space statistics. (b): Free intracellular calcium concentration $[\mathrm{Ca}^{2+}]_\mathrm{i}$ space statistics. (c): Sarcomere length $\mathrm{SL}$ space statistics.}
     \label{fig:spatial_statistics_mmem}
\end{figure}

The results are shown in Figure \ref{fig:results_multiscale_aggregate}, in the form of pressure-volume (PV) loops in Figure \ref{fig:pv_multiscale} and endocardial left ventricular pressure ($p_\mathrm{LV}$) in Figure \ref{fig:p_lv}. 
Figure \ref{fig:pv_multiscale} shows similar behavior of all four models, and it particularly shows almost no appreciable differences between eikonal and monodomain-driven models. The most apparent difference between the four models is that the models incorporating the MCF exhibit higher pressures in the initial and final phases of the heartbeat, when the active force transient has almost completely relaxed. This is due to the higher diastolic active force of the calibrated cellular electromechanical model with respect to the one neglecting the MCF, as shown in Figure \ref{fig:af_calibration}. In Figure \ref{fig:spatial_statistics} we report the active stress, intracellular calcium and sarcomere length comparison for the M and M+MCF models. The active stress and sarcomere lengths of the two models behave similarly, although the diastolic active stress exhibits higher values when incorporating the MCF, as seen in Figure \ref{fig:as_ss}. The intracellular calcium, represented in Figure \ref{fig:intraca_ss}, instead, behaves similarly to the single-cell model, exhibiting a higher peak value, nevertheless maintaining a realistic waveform. In Figure \ref{fig:spatial_statistics_mmem} we report the active stress, intracellular calcium and sarcomere length comparison for the M+MCF and E+MCF models. The differences between the two models is even less noticeable in this case, as the computed active stresses and intracellular calcium values are virtually superimposed, as seen in Figures \ref{fig:as_ss_mmem} and \ref{fig:intraca_ss_mmem}.

\subsection{Effects of perturbation of calcium sensitivity}
\label{sec:numerical_calcium_sens}
In the previous section we have compared four different cardiac electromechanics models, models that were either monodomain- or eikonal-driven and either incorporated the MCF or not. We have seen that in baseline conditions, the four models produce almost indistinguishable PV loops. In this section we examine the effects of the MCF on the electromechanical models in the case of perturbations to the calcium dissociation constant $k_\mathrm{d}$. In Section \ref{sec:cellular_model_perturbed} we examine the effects on the cellular scale, while in Section \ref{sec:multiscale_model_perturbed} we do so for the multiscale model, where we limit ourselves to the eikonal-driven models.
\subsubsection{Cellular model}
\label{sec:cellular_model_perturbed}
\begin{figure}[H]
     \centering
     \subfloat[][]{\includegraphics[width = 0.48\textwidth]{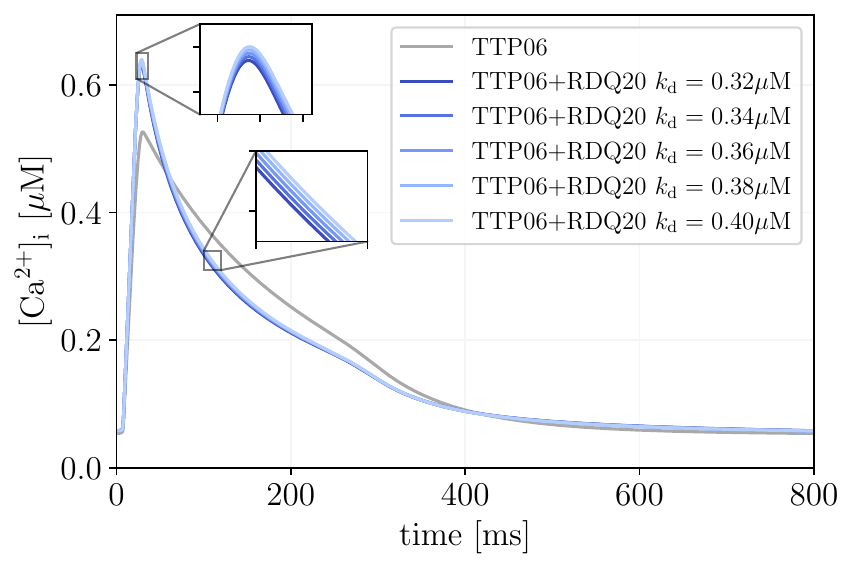}\label{fig:cai_sensitivity}}\hfill
     \subfloat[][]{\includegraphics[width = 0.48\textwidth]{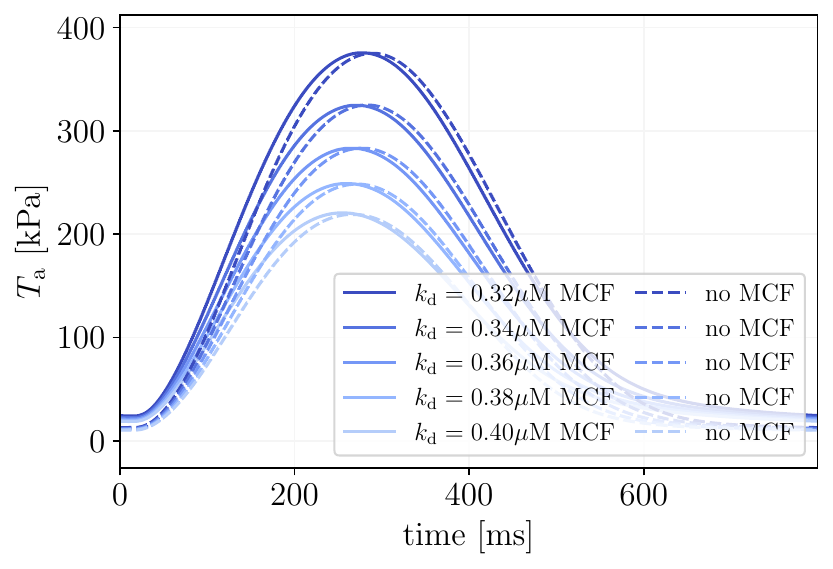}\label{fig:af_sensitivity}}
     \caption{Free intracellular calcium and active force transients depending on the calcium dissociation constant $k_\mathrm{d}$ in the case of the TTP06+RDQ20 model with and without the MCF. Darker blues indicate lower values of $k_\mathrm{d}$, whereas lighter blues indicate higher values. (a): Free intracellular calcium $[\mathrm{Ca}^{2+}]_\mathrm{i}$ depending on $k_\mathrm{d}$. In grey is the TTP06 model without the feedback, which is independent from $k_\mathrm{d}$. Insets represent peak and relaxation. (b): Variation of the active force $T_\mathrm{a}$ depending on $k_\mathrm{d}$. Full lines represent the TTP06+RDQ20 model with the MCF, whereas dashed lines are the TTP06+RDQ20 model without the feedback.}
     \label{fig:sensitivity_cellular}
\end{figure}

\begin{figure}[h]
     \centering
     \subfloat[][]{\includegraphics[width = 0.48\textwidth]{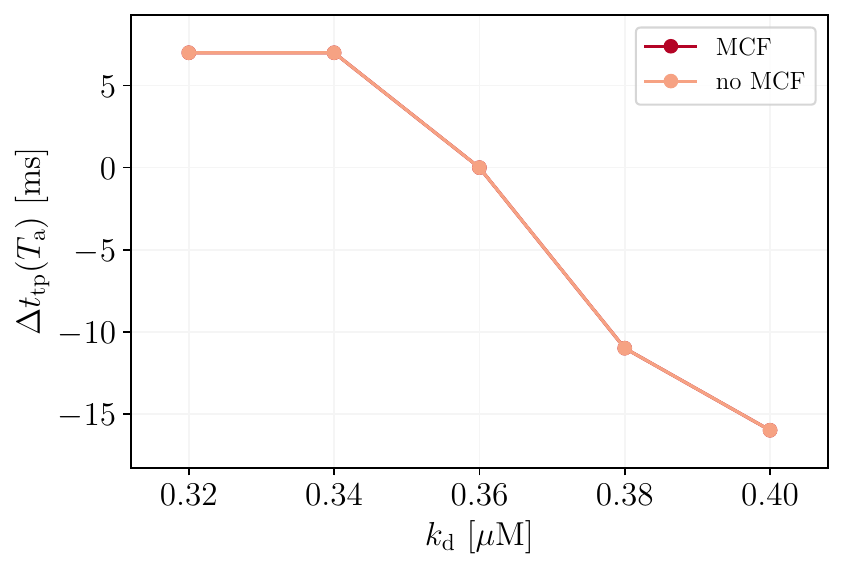}\label{fig:ttps_sensitivity}}\hfill
     \subfloat[][]{\includegraphics[width = 0.48\textwidth]{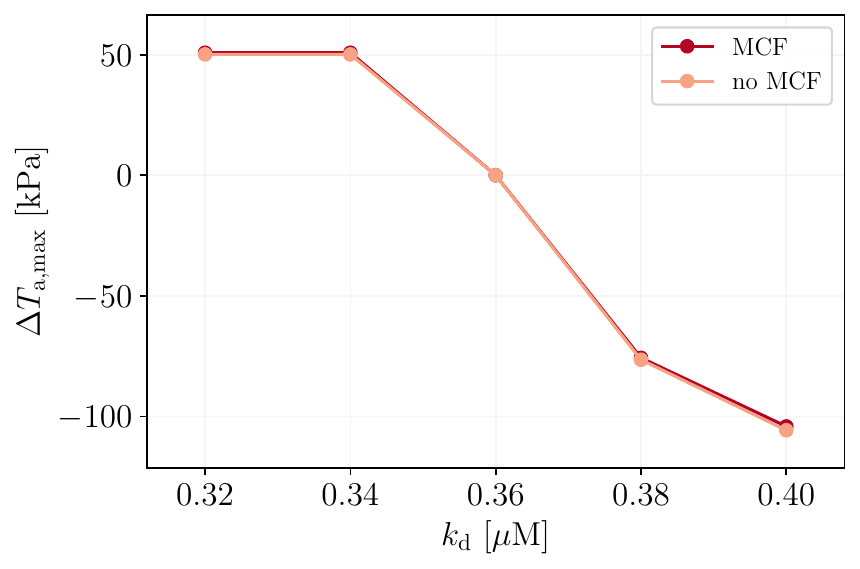}\label{fig:tamax_sensitivity}}
     \caption{Active force transients time to peak and peak value variation depending on the calcium dissociation constant $k_\mathrm{d}$ in the case of the TTP06+RDQ20 model with the MCF (red) and without the MCF (yellow). (a): Variation of the active force time to peak $\Delta t_\mathrm{tp}$ depending on $k_\mathrm{d}$. (b): Variation of the peak active force $\Delta T_\mathrm{a,max}$ depending on $k_\mathrm{d}$.}
     \label{fig:sensitivity_cellular_afkinetics}
\end{figure}

\begin{figure}[h]
     \centering
     \subfloat[][]{\includegraphics[width = 0.48\textwidth]{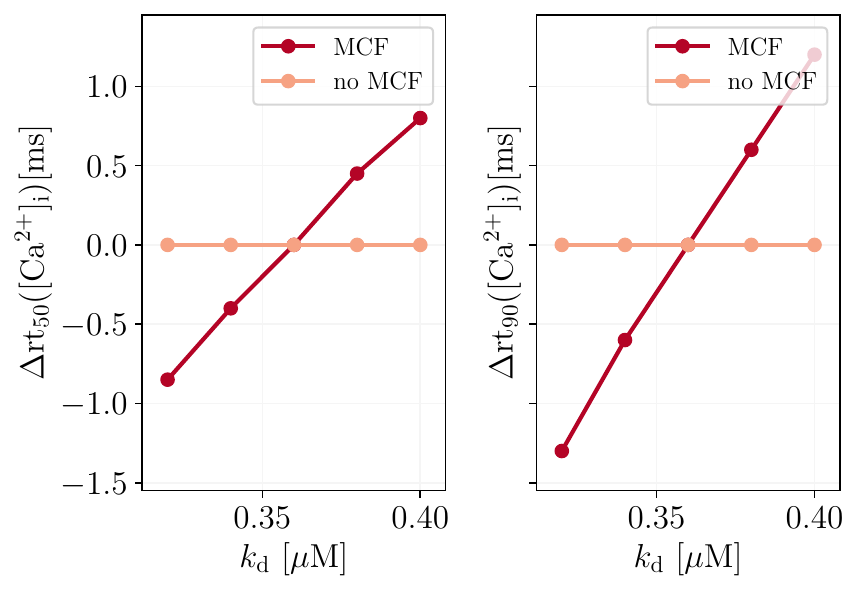}\label{fig:rts_ca}}\hfill
     \subfloat[][]{\includegraphics[width = 0.48\textwidth]{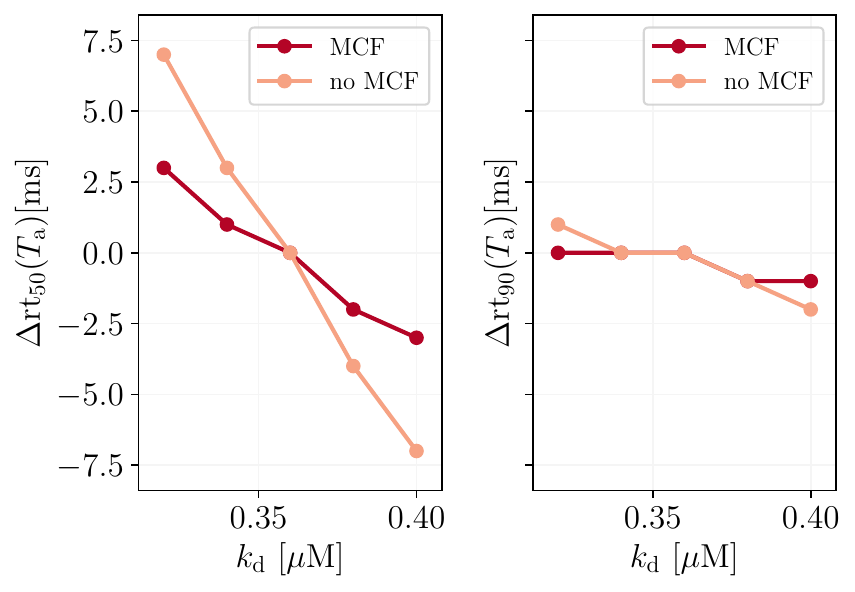}\label{fig:rts_af}}
     \caption{Variations of the relaxation times of the intracellular calcium transient and active force transient depending on the calcium dissociation constant $k_\mathrm{d}$ in the case of the TTP06+RDQ20 model with the MCF (red) and without the MCF (yellow). (a): Variation of the intracellular calcium amplitude $50\%$ and $90\%$ relaxation times, $\mathrm{rt}_{50}([\mathrm{Ca}^{2+}]_\mathrm{i})$ and $\mathrm{rt}_{90}([\mathrm{Ca}^{2+}]_\mathrm{i})$, depending on $k_\mathrm{d}$. (b): Variation of the active force amplitude $50\%$ and $90\%$ relaxation times, $\mathrm{rt}_{50}(T_\mathrm{a})$ and $\mathrm{rt}_{90}(T_\mathrm{a})$, depending on $k_\mathrm{d}$.}
     \label{fig:relaxation_times_sens}
\end{figure}

The results for the perturbed cellular models are reported in Figure \ref{fig:sensitivity_cellular} in the form of free intracellular calcium concentration transients and microscopic active force transients. 
As a first consequence of the feedback \eqref{eq:catn_dt_active_force}, we have that the perturbation of the calcium dissociation constant $k_\mathrm{d}$ affects the intracellular calcium transient, reported in Figure \ref{fig:cai_sensitivity}. Indeed, for higher values of the calcium dissociation constant $k_\mathrm{d}$, as expected from Equation \eqref{eq:transition_bound}, the equilibrium for the calcium-troponin binding reaction shifts to the left toward an unbound state for the calcium, raising the overall free intracellular calcium concentration. This effect cannot be seen in the case where no feedback is modeled. Regarding the active force, Figure \ref{fig:af_sensitivity} shows that for higher $k_\mathrm{d}$ values, the developed active force is lower, both for the case with and without MCF, due to the lower propensity of the calcium to be bound to troponin, and consequently lower permissivity of the regulatory units, as expected \cite{regazzoni_biophysically_2020}. The magnitude of this effect is similar between the two cases, as exemplified by Figure \ref{fig:sensitivity_cellular_afkinetics}, where we compare the variations under the combined effects of the MCF and perturbation of parameter $k_\mathrm{d}$ of time to peak of the active force $\Delta t_\mathrm{tp}(T_\mathrm{a})$ and of the peak active force $\Delta T_\mathrm{a,max}$ with respect to the baseline value with the MCF, computed as:
\begin{equation}
    \Delta t_\mathrm{tp}(T_\mathrm{a};k_\mathrm{d}) = t_\mathrm{tp}(T_\mathrm{a};k_\mathrm{d}) - t_\mathrm{tp}(T_\mathrm{a};k_{\mathrm{d},0}), \qquad \Delta T_\mathrm{a,max}(k_\mathrm{d}) = T_\mathrm{a,max}(k_\mathrm{d}) - T_\mathrm{a,max}(k_{\mathrm{d},0}),
\end{equation}
where $k_{\mathrm{d},0}=0.36\si{\micro\Molar}$. From Figures \ref{fig:ttps_sensitivity} and \ref{fig:tamax_sensitivity}, we see again the small shift in the wavefront in the case of the presence of the feedback, as already seen in the calibrated model in Section \ref{sec:cellular_electromechanical_model}. This shift remains consistent across variations of the calcium dissociation constant $k_\mathrm{d}$ for the cellular model, and has virtually no impact on the time to peak and maximum active force trend. Regarding the relaxation times of the intracellular calcium transient and active force, we measure the combined effect by computing their variations $\Delta\mathrm{rt}_{50/90}([\mathrm{Ca}^{2+}]_\mathrm{i})$ and $\Delta\mathrm{rt}_{50/90}(T_\mathrm{a})$, as:
\begin{equation}
\begin{aligned}
    \Delta\mathrm{rt}_{50/90}([\mathrm{Ca}^{2+}]_\mathrm{i};k_{\mathrm{d}}) &= \mathrm{rt}_{50/90}([\mathrm{Ca}^{2+}]_\mathrm{i};k_\mathrm{d}) - \mathrm{rt}_{50/90}([\mathrm{Ca}^{2+}]_\mathrm{i};k_{\mathrm{d},0}), \\ \Delta\mathrm{rt}_{50/90}(T_\mathrm{a};k_{\mathrm{d}}) &= \mathrm{rt}_{50/90}(T_\mathrm{a};k_\mathrm{d}) - \mathrm{rt}_{50/90}(T_\mathrm{a};k_{\mathrm{d},0}),
\end{aligned}
\end{equation}
where again $k_{\mathrm{d},0}=0.36\si{\micro\Molar}$. The variations are reported in Figure \ref{fig:relaxation_times_sens}, where Figure \ref{fig:rts_ca} reflects the fact that the perturbation of the calcium dissociation constant has no effect when the MCF is not included, whereas its increase prolongs the relaxation times. Figure \ref{fig:rts_af} instead shows that the inclusion of the MCF causes the variations of the relaxation times of the active force transient to be less prominent in isometric conditions.

\subsubsection{Multiscale model}
\label{sec:multiscale_model_perturbed}
\begin{figure}[ht]
     \centering
     \subfloat[][]{\includegraphics[width = 0.48\textwidth]{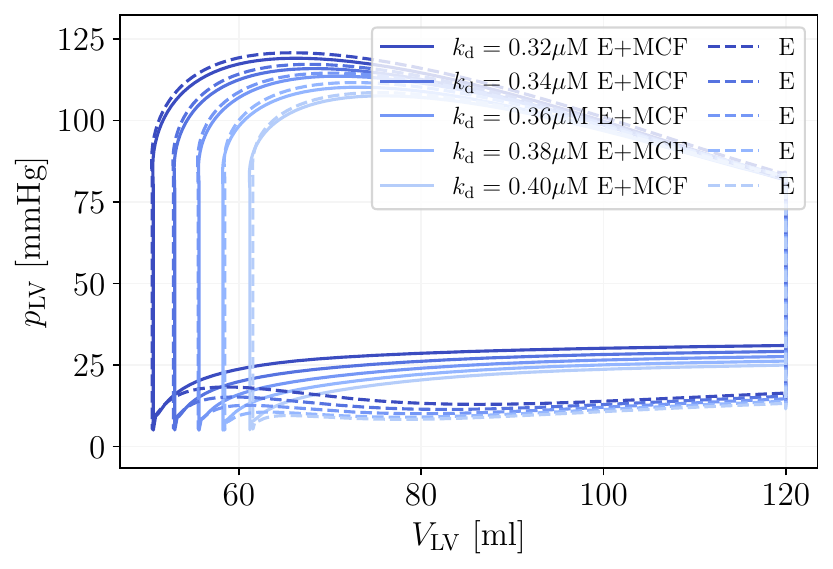}\label{fig:pvloops_sensitivity}}\hfill
     \subfloat[][]{\includegraphics[width = 0.48\textwidth]{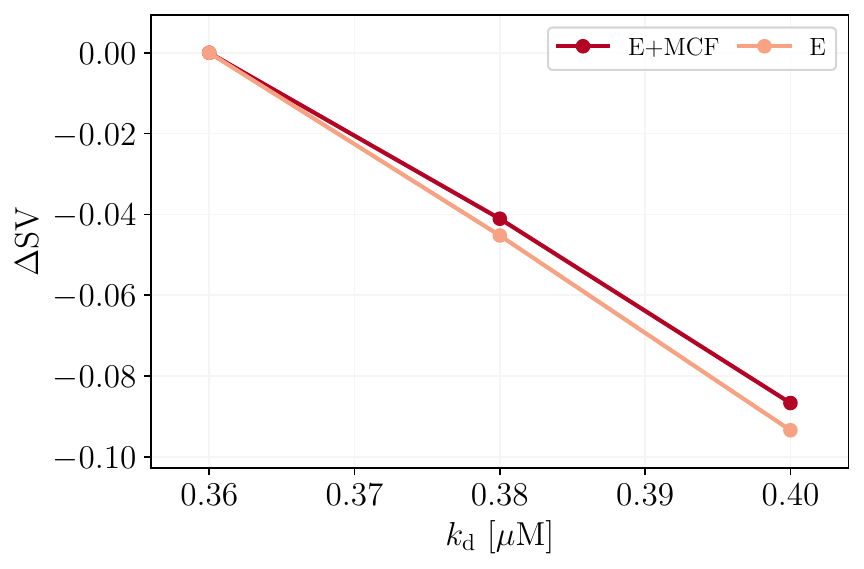}\label{fig:metrics_sensitivity}}
     \caption{Multiscale results in terms of PV loops and depending on the calcium dissociation constant $k_\mathrm{d}$. (a): PV loops for varying $k_\mathrm{d}$, where darker blues indicate lower $k_\mathrm{d}$ values. Full lines represent the TTP06+RDQ20 model with the MCF, whereas dashed lines are the TTP06+RDQ20 model without the feedback. (b): Variation of $\Delta\mathrm{SV}$ dor increasing values of $k_\mathrm{d}$, for both TTP06+RDQ20 model with the MCF (red) and without the MCF (yellow).}
     \label{fig:sensitivity_multiscale}
\end{figure}

Regarding the multiscale model, the responses to the perturbation of $k_\mathrm{d}$ are reported in Figure \ref{fig:sensitivity_multiscale}. For the two models the response to parameter perturbation have the same effects although of different magnitudes. In order to better quantify the differences between the effects for both cases, we compute the relative variations of the stroke volume (SV) with respect to the baseline case as:
\begin{equation}
    \Delta \mathrm{SV}(k_\mathrm{d}) = \frac{\mathrm{SV}(k_\mathrm{d}) - \mathrm{SV}(k_{\mathrm{d},0})}{\mathrm{SV}(k_{\mathrm{d},0})},
\end{equation}
where $k_\mathrm{d,0}=0.36\si{\micro \Molar}$. We report $\Delta\mathrm{SV}$ as as a function of $k_\mathrm{d}$ in Figure \ref{fig:metrics_sensitivity}. We can observe how the stroke volume decreases for increasing $k_\mathrm{d}$, but this effect is less pronounced if MCF is included. We recall, from Section \ref{sec:cellular_model_perturbed}, that in single-cell isometric simulations the perturbation of $k_\mathrm{d}$ produces very similar responses in terms of maximal active force and time to peak, irrespective of the inclusion of the MCF. Instead, the shortening of the relaxation times at higher $k_\mathrm{d}$ values is less evident when including the MCF. This makes it so that the observed differences emerging in the multiscale model could be owed to the overall longer active force transient durations as well as the heterogeneity and dynamic properties of the multiscale problem combined with the MCF.
\subsection{Remark on the computational times}
\begin{figure}[ht]
     \centering
     \subfloat[][]{\includegraphics[width = 0.48\textwidth]{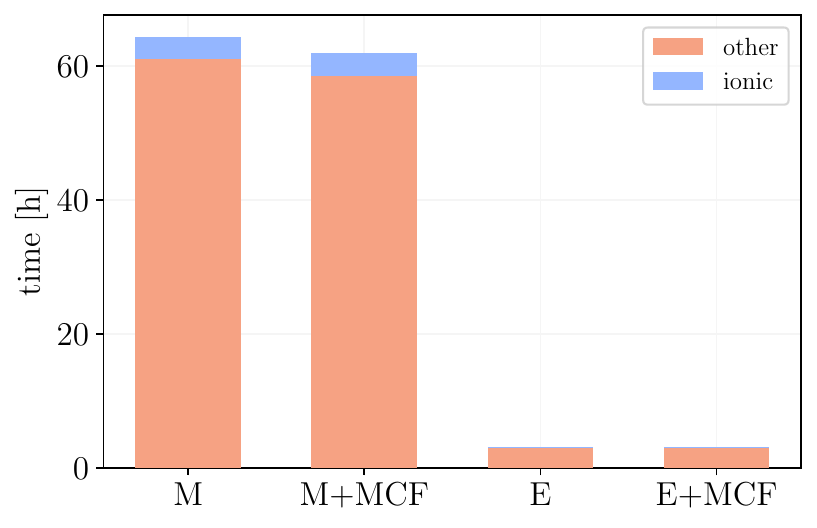}\label{fig:comptimes}}\hfill
     \subfloat[][]{\includegraphics[width = 0.48\textwidth]{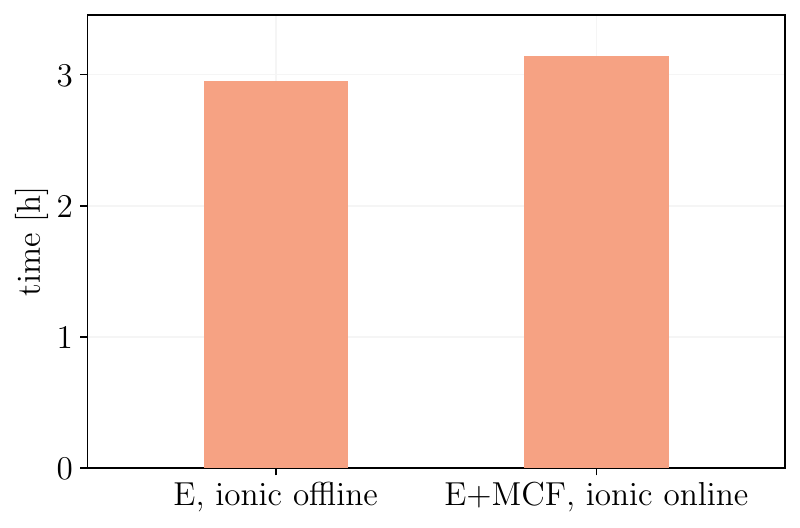}\label{fig:comptimes_wc}}
     \caption{Allocation of computational times for all the multiscale models in Section \ref{sec:numerical_results_multiscale}, on 72 cores. (a): Total computational times on 72 cores for the multiscale models. In blue is the time for the solution of the ionic models, in red for all the other operations. Total bar height is the total computational time. (b): Computational times of eikonal-driven models, in the case of no feedback without considering the ionic model solution times (left) and in the case with the MCF and online ionic solution (right).}
     \label{fig:all_comptimes}
\end{figure}

All computational times are reported in Figure \ref{fig:all_comptimes}. The simulation of two heartbeats in Section \ref{sec:numerical_results_multiscale} of the eikonal-driven simulations, E and E+MCF, both took $3\si{\hour}\ 8\si{\minute}$ to simulate on 72 cores, where in both cases we solved the ionic models online. For both cases the solution of the ionic model took $11\si{\minute}$ (corresponding to about 6\% of the total time), which we take to be the most conservative estimate of the additional computational time of the online ionic solution with respect to the offline one.
\par
Conversely, the monodomain-driven simulations, M and M+MCF, took $61\si{\hour}\ 56\si{\minute}$ for M and $64\si{\hour}\ 26\si{\minute}$ for M+MCF on 72 cores, not including the time used for output and the construction times for the RBF interpolator for the common quantities ($\mathbf{F}$, $\mathbf{d}$, and $w_\mathrm{Ca}$), implying no particular additional computational burden for the simulations including the MCF.
\par
It is evident that eikonal-driven models offer a significant speedup on simulation runtime with respect to their monodomain-driven counterparts, as highlighted in Figure \ref{fig:comptimes}. Moreover, for the eikonal-driven simulations the additional overhead due to the necessary solution of the ionic models online is low, even in the most conservative scenario, as reported in Figure \ref{fig:comptimes_wc}.
\section{Discussion and conclusions}
\label{sec:discussion}
{In this work we have investigated the impact of the inclusion of the MCF in multiscale cardiac electromechanics simulations on the overall cardiac function. To achieve this, we have derived a fully coupled cellular electromechanical model by coupling a ionic model and a force generation model at the subcellular scale. We extended this derivation to a broader class of ionic models, including single-buffer models, capable of retaining realistic cellular behavior. Specifically, we used this procedure to obtain a cellular electromechanical model starting from the TTP06 model. The cellular models revealed an increased influence of ionic and activation parameters on calcium and active force dynamics, due to the bidirectional coupling. Due to this increased sensitivity, a naive coupling without parameter calibration may lead to inconclusive results about the role of the MCF in multiscale simulations. To resolve this issue, we presented a calibration procedure suitable for calibrating the model's parameters in order to obtain equivalent baseline conditions at the cellular scale. Such a procedure is able to recover realistic calcium transients and active force kinetics, which are important metrics when evaluating an electromechanical model's biophysical fidelity. } 
\par
{Starting from the calibrated bidirectionally coupled cellular electromechanical model, we investigated its multiscale extension to a fully coupled electromechanical model, capable of reproducing the MCF. We have seen that the inclusion of the MCF in baseline conditions has little effect on overall cardiac function, although there is no evident additional computational cost of including the MCF with respect to standard, geometrically-coupled multiscale electromechanical simulations.} As mentioned in Section \ref{sec:introduction}, the solution of problems arising from cardiac electromechanics, such as Equations \eqref{eq:full_em_space}, suffers from several computational bottlenecks, which are amplified by the coupled nature of the system. Equations \eqref{eq:ionic_variables_space_mcf} and \eqref{eq:monodomain} require high spatial resolutions and fast dynamics, imposing small timesteps and fine meshes in order to capture sharp depolarization wavefronts. Consequently, the high number of degrees of freedom is carried onto the approximation of Equations \eqref{eq:mechanics} and \eqref{eq:sarcomere_variables_space}. Even when intergrid methods are used in order to alleviate these problems \cite{salvador_intergrid_2020,bucelli_preserving_2023,bucelli_robust_2024}, the computational burden of these is non-negligible, in terms of memory and computational time. This is why it is important to find modeling frameworks which intrinsically require lower computational resources while guaranteeing a high degree of biophysical fidelity, such as the one proposed in Section \ref{sec:eikonal_driven}, where in alternative to the fully coupled multiscale framework, we proposed an eikonal-driven framework. We have demonstrated the eikonal-driven model to be equally appropriate for capturing the MCF and subsequent aggregate electromechanical behavior, similarly to eikonal-based electromechanical simulations neglecting the MCF \cite{stella_fast_2022}. Moreover, we have shown that the eikonal-driven model including the MCF behaves differently from the models neglecting the feedback also in the eikonal-driven case, justifying its use in medically relevant scenarios, such as inotropic modulation.
\par
Due to the formulation of the eikonal model employed in the present study, the present results are limited to sinus rhythm. However, recent advancements in eikonal-based modeling have introduced formulations capable of modeling reentries \cite{barrios_espinosa_cyclical_2025, gander_eikonal_2024, jacquemet_eikonal_2010}, for which the inclusion of the MCF can be extended. Additionally, since the calibration of the cellular model relies on a grid search algorithm, the precision of the calibration is limited by the refinement level of the parameter grids, The accuracy of the cellular electrophysiology is inherited from the TTP06 model. However, the TTP06 model is still a widely used tool in multiscale cardiac electromechanics simulations, and the presented analysis can be used as a starting point for more complex ionic models.
\section*{Acknowledgments}
\label{sec:acknowledgemnts}
The present research has been partially supported by MUR, grant Dipartimento di Eccellenza 2023-2027.
F.R. and S.P. have received support from the project PRIN2022, MUR, Italy, 2023-2025, P2022N5ZNP “SIDDMs: shape-informed data-driven models for parametrized PDEs, with application to computational cardiology”, funded by the European Union (Next Generation EU, Mission 4 Component 2). M.B. and L.D. have received support from the project  PRIN2022, MUR, Italy, 2023-2025, proj. no. 202232A8AN ``Computational modeling of the human heart: from efficient numerical solvers to cardiac digital twins". I.R., F.R., M.B., S.P., L.D. are members of GNCS, “Gruppo Nazionale per il Calcolo Scientifico" (National Group for Scientific Computing) of INdAM (Istituto Nazionale di Alta Matematica). F.R., M.B., S.P., L.D. have received support from the EuroHPC JU project dealii-X (grant number 101172493) funded under the HORIZON-EUROHPC-JU-2023-COE-03-01 initiative. I.R., F.R., M.B., S.P., L.D. acknowledge the INdAM GNCS project CUP E53C24001950001. The authors acknowledge ISCRA for awarding this project access to the LEONARDO supercomputer, owned by the EuroHPC Joint Undertaking, hosted by CINECA (Italy).
\appendix
\section{Model and numerical parameters}
\label{app:model_numerical_parameters}
We report the model and numerical parameters used in the simulations, if different from what is reported in previous works \cite{fedele_comprehensive_2023,regazzoni_cardiac_2022,regazzoni_machine_2020}.
\begin{table}[H]
    \centering 
    \begin{tabular}{ c | c | c | c }
    Parameter & Unit & TTP06+RDQ20 & RDQ20 \\ 
    \hline\hline
    $\mu$ & $-$ & 10 & 10 \\
    $\gamma$ & $-$ & 30 & 30 \\ 
    $Q$ & $-$ & 2 & 2 \\ 
    ${k}_\mathrm{d}$ & $[\si{\micro\Molar}]$ & 0.36 & 0.36\\ 
    $\alpha_{{k}_\mathrm{d}}$ & $\si[per-mode=reciprocal]{\micro\Molar \per \micro \meter}$ & -0.2083 & -0.2083 \\ 
    ${k}_\mathrm{off}$ & $\si[per-mode=reciprocal]{\second}$ & 5 & 8\\
    ${k}_\mathrm{basic}$ & $\si[per-mode=reciprocal]{\second}$ & 5 & 4 \\
    $\mu_{f_\mathcal{P}}^0$ & $\si[per-mode=reciprocal]{\second}$ & 32.255 & 32.255 \\
    $\mu_{f_\mathcal{P}}^1$ & $\si[per-mode=reciprocal]{\second}$ & 0.768 & 0.768 \\
    $r_0$ & $\si[per-mode=reciprocal]{\second}$ & 134.31 & 134.31  \\
    $\alpha$ & $-$ & 25.184 & 25.184 \\
    $a_\mathrm{XB}$ & $\si[per-mode=reciprocal]{\mega \pascal}$ & 3319.65 & 2250 \\
    \end{tabular}
    \\[10pt]
    \caption{RDQ20 model parameters in the case of the calibrated TTP06+RDQ20 model and the baseline simulation (denoted here only with RDQ20).}
    \label{table:rdq20_parameters_calibrations_appendix}
\end{table}

\begin{table}[H]
    \centering 
    \begin{tabular}{ c | c | c | c }
    Parameter & Unit & Description & Value \\ 
    \hline\hline
    $p_{AVO}$ & $\si{\pascal}$ & Aortic valve opening pressure & 11000\\
    $p_{MVO}$ & $\si{\pascal}$ & Mitral valve opening pressure & 667\\ 
    $C$ & $\si[per-mode=reciprocal]{\cubic\meter\per \pascal}$ & Vessel capacitance & $\num[exponent-product=\ensuremath{\cdot}]{4.5e-9}$\\ 
    $R$ & $\si[per-mode=reciprocal]{\pascal\second\per\cubic\meter}$ & Distal resistance & $\num[exponent-product=\ensuremath{\cdot}]{5.5e7}$\\
    $p_{ED}$ & $\si[per-mode=reciprocal]{\pascal}$ & End-diastolic pressure & $\num[exponent-product=\ensuremath{\cdot}]{1333.0}$ \\ 
    $V_{ED}$ & $\si[per-mode=reciprocal]{\milli\litre}$ & End-diastolic volume & $\num[exponent-product=\ensuremath{\cdot}]{120}$\\
    \end{tabular}
    \\[10pt]
    \caption{Circulation and Windkessel model parameters.}
    \label{table:circulation_parameters_calibrations_appendix}
\end{table}

\begin{table}[H]
    \centering 
    \begin{tabular}{ c | c | c | c }
    Parameter & Unit & Description & Value \\ 
    \hline\hline
    $\sigma_f$ & $\si[per-mode=reciprocal]{\meter\squared\per\second}$ & Longitudinal conductivity & \num[exponent-product=\ensuremath{\cdot}]{1.6603e-4} \\
    $\sigma_s$ & $\si[per-mode=reciprocal]{\meter\squared\per\second}$ & Transversal conductivity & \num[exponent-product=\ensuremath{\cdot}]{0.7590e-4} \\
    $\sigma_n$ & $\si[per-mode=reciprocal]{\meter\squared\per\second}$ & Normal conductivity & \num[exponent-product=\ensuremath{\cdot}]{0.2443e-4}\\
    \hline
     \begin{tabular}{@{}c@{}} $(x_0,y_0,z_0)$ \\ $(x_1,y_1,z_1)$ \\ $(x_2,y_2,z_2)$\end{tabular} & $\si[per-mode=reciprocal]{\meter}$ & Impulse sites (M+E) &\begin{tabular}{@{}c@{}} $(0.04229, 1.34726, 0.05256)$ \\ $(0.07065, 1.36207, 0.04283)$ \\ $(0.06783, 1.31976, 0.04419)$\end{tabular}\\ 
     $\mathcal{I}_\mathrm{app}$ & $\si{\second}$ & Impulse duration (M) & $\num[exponent-product=\ensuremath{\cdot}]{3e-3}$\\
     $\mathcal{I}_\mathrm{app}$ & $\si[per-mode=reciprocal]{\volt\per\second}$ & Impulse amplitude (M) &  34.28 \\ $\mathcal{I}_\mathrm{app}$ & $\si{\meter}$ & Spherical impulse radius (M) &  $\num[exponent-product=\ensuremath{\cdot}]{2.5e-3}$ \\
     $t_\mathrm{app}$ & $\si{\second}$ & Depolarization current duration (E) &  $\num[exponent-product=\ensuremath{\cdot}]{2e-3}$ \\
     $\Bar{\mathcal{I}}_\mathrm{app}$ & $\si[per-mode=reciprocal]{\volt\per\second}$ & Depolarization current amplitude (E) &  $\num[exponent-product=\ensuremath{\cdot}]{25.71}$ \\
     $c_0$ & $\si[per-mode=reciprocal]{{\second}^{-1/2}}$ & Wavefront velocity parameter (E) &  $\num[exponent-product=\ensuremath{\cdot}]{52.195}$ \\
    \end{tabular}
    \\[10pt]
    \caption{Monodomain (M) and eikonal/ionic (E) model parameters.}
    \label{table:monodomain_appendix}
\end{table}

\begin{table}[H]
    \centering 
    \begin{tabular}{ c | c | c | c }
    Parameter & Unit & Description & Value \\ 
    \hline\hline
    $\rho_s$ & $\si[per-mode=reciprocal]{\kilogram\per\cubic\meter}$ & Solid density & \num[exponent-product=\ensuremath{\cdot}]{1000} \\
    $K_\perp$ & $\si[per-mode=reciprocal]{\pascal\per\meter}$ & Normal stiffness on $\Gamma^\mathrm{epi}$  & 200000\\
    $K_\parallel$ & $\si[per-mode=reciprocal]{\pascal\per\meter}$ & Tangential stiffness on $\Gamma^\mathrm{epi}$ & 20000\\
    $C_\perp$ & $\si[per-mode=reciprocal]{\pascal\second\per\meter}$ & Normal viscosity on $\Gamma^\mathrm{epi}$  & 20000 \\
    $C_\parallel$ & $\si[per-mode=reciprocal]{\pascal\second\per\meter}$ & Tangential viscosity on $\Gamma^\mathrm{epi}$  & 2000\\
    \end{tabular}
    \\[10pt]
    \caption{Mechanics model parameters.}
    \label{table:mechanics_appendix}
\end{table}

% TODO passive mechanics active stress electrophysiology
\begin{table}[H]
    \centering 
    \begin{tabular}{ p{18em} | c | c }
    Parameter & M/M+MCF & E/E+MCF \\
    \hline\hline
    Element type & Tet & Tet \\
    $\mathcal{T}_h$ average cell diameter & $0.67 \si[per-mode=reciprocal]{\milli\meter}$ & $1.15 \si[per-mode=reciprocal]{\milli\meter}$ \\
    $\mathcal{T}_H$ average cell diameter & $1.15 \si[per-mode=reciprocal]{\milli\meter}$ & - \\
    Electrophysiology FE degree & 2 & 1  \\
    Active stress FE degree & 1 & 1  \\
    Mechanics FE degree & 1 & 1 \\
    Electrophysiology timestep $\tau$ & $\num[exponent-product=\ensuremath{\cdot}]{5e-5}\si{\second}$ & $\num[exponent-product=\ensuremath{\cdot}]{5e-5}\si{\second}$  \\
    Active stress and Mechanics timestep $\Delta t$ & $\num[exponent-product=\ensuremath{\cdot}]{0.001}\si{\second}$ & $\num[exponent-product=\ensuremath{\cdot}]{0.001}\si{\second}$  \\
    Heartbeat duration $T$ & $\num[exponent-product=\ensuremath{\cdot}]{0.8}\si{\second}$ & $\num[exponent-product=\ensuremath{\cdot}]{0.8}\si{\second}$  \\
    Final time $T_\mathrm{fin}$ & $\num[exponent-product=\ensuremath{\cdot}]{1.6}\si{\second}$ & $\num[exponent-product=\ensuremath{\cdot}]{1.6}\si{\second}$  \\
    \end{tabular}
    \\[10pt]
    \caption{Numerical discretization parameters for the multiscale electromechanics problem.}
    \label{table:numerical_discretization_appendix}
\end{table}

\begin{table}[H]
    \centering 
    \begin{tabular}{ p{12em} | c | c }
    Interpolated quantity & $M$ & $\alpha$ \\
    \hline\hline
    $w_\mathrm{Ca}$ & $6$ & $1.75$\\
    $\mathbf{d}$ & $6$ & $1.75$ \\
    $\mathbf{F}$ (SVD) & $4$ & $1.75$ \\
    $y_{\mathrm{dCa_{Tn}}}$ & $6$ & $1.75$ \\
    \end{tabular}
    \\[10pt]
    \caption{Parameters used for the RBF interpolation, when used, reported as in \cite{bucelli_preserving_2023}.}
    \label{table:numerical_discretization_rbf}
\end{table}

\bibliographystyle{elsarticle-num} 
\bibliography{cas-refs}



%% For numbered reference style

\end{document}